\documentclass[10pt]{article}
\usepackage{sectsty,latexsym,amssymb,dsfont,textcomp}
\usepackage[centertags]{amsmath}
\usepackage{fullpage}
\usepackage[centertags]{amsmath}
\sectionfont{\normalsize}
\subsectionfont{\normalsize}
\begin{document}
\date{}
\numberwithin{equation}{section}
\title{Weak-strong uniqueness of dissipative measure-valued solutions 
for polyconvex elastodynamics
}

\author{
Sophia Demoulini\footnote{Centre for Mathematical Sciences, University of
Cambridge, UK, S.Demoulini@dpmms.cam.ac.uk} 
\and
David M.A. Stuart\footnote{Centre for Mathematical Sciences, University of
Cambridge, UK, D.M.A.Stuart@damtp.cam.ac.uk} 
\and 
Athanasios E. Tzavaras\footnote{Department of Applied Mathematics, University of 
Crete and Institute for Applied and Computational Mathematics, FORTH\@, Heraklion, Greece,
tzavaras@tem.uoc.gr}
}

\maketitle


%
%

\newcommand{\wkarr}{\; \rightharpoonup \;}
\def\Weak{\,\,\relbar\joinrel\rightharpoonup\,\,}
\def\del{\partial}
\font\msym=msbm10
\def\Real{{\mathop{\hbox{\msym \char '122}}}}
\def\R{\Real}
\def\A{\mathbb A}
\def\Z{\mathbb Z}
\def\K{\mathbb K}
\def\J{\mathbb J}
\def\L{\mathbb L}
\def\D{\mathbb D}
\def\Mink{{\mathop{\hbox{\msym \char '115}}}}
\def\Integers{{\mathop{\hbox{\msym \char '132}}}}
\def\Complex{{\mathop{\hbox{\msym\char'103}}}}
\def\C{\Complex}
\font\smallmsym=msbm7
\newcommand{\beq}{\begin{equation}}
\newcommand{\eeq}{\end{equation}}
\newcommand\la{\label}
\newcommand{\id}{\mathds{1}}

\newtheorem{lemma}{Lemma}[section]
\newtheorem{theorem}[lemma]{Theorem}
\newtheorem{maintheorem}[lemma]{Main Theorem}
\newtheorem{corollary}[lemma]{Corollary}
\newtheorem{definition}[lemma]{Definition}
\newtheorem{remark}[lemma]{Remark}
\newtheorem{remarks}[lemma]{Remarks}
\newtheorem{Notation}[lemma]{Notation}
\newcommand{\proof}{\noindent {\it Proof}\;\;\;}
\newcommand{\qed}{\protect~\protect\hfill $\Box$}

\newcommand{\cd}{{\cal D}}
\newcommand{\sma}{_{{ A}}}
\newcommand{\ce}{{\cal E}}
\newcommand{\dist}{{\mbox{dist}}}
\newcommand{\cs}{{\cal S}}
\newcommand{\ca}{{\cal A}}
\newcommand{\cc}{{\cal C}}
\newcommand{\cp}{{\cal P}}
\newcommand{\cb}{{\cal B}}
\newcommand{\cl}{{\cal L}}
\newcommand{\cg}{{\cal G}}
\newcommand{\cq}{{\mathcal Q}}
\newcommand{\hj}{^{J}}
\newcommand{\n}{\mbox{\boldmath $ \nu$}}
\newcommand{\kap}{\mbox{\boldmath $\kappa$}}
\newcommand{\s}{\mbox{\boldmath $ \sigma$}}
\newcommand{\g}{\mbox{\boldmath $ \gamma$}}
\newcommand{\zet}{\mbox{\boldmath $ \zeta$}}
\newcommand{\m}{\mbox{\boldmath $ \mu$}}
\newcommand{\ebb}{\mbox{\boldmath $\epsilon$}}
\newcommand{\nutx}{\n_{t,x}}
\newcommand{\mutx}{\m_{t,x}}
\newcommand{\hjm}{^{J-1}}
\newcommand{\cu}{{\cal U}}
\newcommand{\tn}{{\tilde\|}}
\newcommand{\rbar}{\overline{r}}
\newcommand{\dxdt}{\; dxdt}
\newcommand{\dx}{\; dx}
\newcommand{\oeps}{\overline{\varepsilon}}
\newcommand{\cgl}{\hbox{Lie\,}{\cal G}}
\newcommand{\ttd}{{\tt d}}
\newcommand{\ttdel}{{\tt \delta}}
\newcommand{\ns}{\nabla_*}
\newcommand{\cM}{{\cal M}}
\newcommand{\spsi}{_{{{\Psi}}}}
\newcommand{\kr}{\hbox{Ker}}
\newcommand{\qinfo}{\stackrel{\circ}{Q}_{\infty}}
\newcommand{\be}{\begin{equation}}
\newcommand{\ee}{\end{equation}}
\newcommand{\ba}{\begin{eqnarray}}
\newcommand{\bes}{\[}
\newcommand{\ees}{\]}
\newcommand{\bas}{\begin{eqnarray*}}
\newcommand{\eas}{\end{eqnarray*}}
\newcommand{\ou}{\overline u}
\newcommand{\ov}{\overline v}
\newcommand{\proj}{{P}}
\newcommand{\mat}{\hbox{Mat}}
\newcommand{\tp}{{\tilde P}}
\newcommand{\cof}{\hbox{cof}\,}
\newcommand{\dt}{\hbox{det}\,}
\newcommand{\p}{\partial}
\newcommand{\pt}{\partial_t}
\newcommand{\vh}{{\hat v}}
\newcommand{\xh}{\widehat{\Xi}}
\newcommand{\pa}{\partial_{\alpha}}
\newcommand{\px}{\partial_{x}}
\newcommand{\fh}{\widehat{F}}
\newcommand{\e}{\epsilon}
\newcommand{\V}{\mbox{Var}}
\newcommand{\lra}{\longrightarrow}
\newcommand{\lav}{\lambda_v}
\newcommand{\lavi}{\lambda_{v_i}}
\newcommand{\laX}{\lambda_\Xi}
\newcommand{\laXB}{\lambda_{\Xi^B}}
\newcommand{\laF}{\lambda_F}
\newcommand{\laZ}{\lambda_Z}
\newcommand{\law}{\lambda_w}
\newcommand{\barQT}{{\overline Q}_T}

\newcommand{\eps}{\varepsilon}
\newcommand{\poseul}{x}
\newcommand{\poslag}{X}
\newcommand{\wxioe}{W(\frac{\poseul_{i+1} - \poseul_i}{\e})}
\newcommand{\wpxioe}{W'(\frac{\poseul_{i+1} - \poseul_i}{\e})
-W'(\frac{\poseul_i - \poseul_{i-1}}{\e})}
\newcommand{\sig}{\sum_{i=0}^{N-1}}
\newcommand{\dotxis}{\dot{X}_i^2}\
\newcommand{\dtye}{\dot{\tilde{y}}^\e}
\newcommand{\veleul}{v}
\renewcommand{\theequation}{\thesection.\arabic{equation}}
\font\msym=msbm10
\def\Real{{\mathop{\hbox{\msym \char '122}}}}
\def\Integers{{\mathop{\hbox{\msym \char '132}}}}
\font\smallmsym=msbm7
\def\smr{{\mathop{\hbox{\smallmsym \char '122}}}}
\def\torus{{{\text{\rm T}} \kern-.42em {\text{\rm T}}}}
\def\T3{\torus^3}
\def\div{\hbox{div}\,}
%
%

%

\begin{abstract}
\noindent
For the equations of elastodynamics with polyconvex stored energy, 
and some related simpler systems,
we define a notion of dissipative measure-valued solution 
and show that such a solution 
agrees with a classical solution with the same initial data 
when such a classical solution exists.  As an application of the
method we give a short proof of strong convergence in the continuum limit
of a lattice approximation of one dimensional elastodynamics
in the presence of a classical solution.  
Also, for a system of conservation laws  endowed with a positive and 
convex entropy, we show 
that dissipative measure-valued solutions attain their initial data 
in a strong sense after time averaging.
\end{abstract}
%
%
%
\section{Introduction}
\label{intro}
\setcounter{equation}{0}
In this article we consider the system of equations of elastodynamics, 
with a stored energy function
which satisfies the condition of polyconvexity introduced in \cite{ball77}. This system
can be embedded into a symmetrizable hyperbolic system of conservation laws
which admits a convex entropy(\cite{dst2, daf, qin}). Using this embedding, the 
existence of globally defined measure-valued solutions (that satisfy
additional geometric properties involving the null Lagrangians) was 
proved in \cite{dst2}, using a method of variational approximation.
The concept of measure-valued solution was introduced into the theory of 
conservation laws in \cite{diperna},
and then into the theory of the incompressible Euler equations in \cite{dipm},
after the 
development of Young measures and weak convergence methods for partial differential equations (\cite{tartar79,evans}).
For several equations of mathematical physics it is currently the only 
notion of solution which is sufficiently broad to allow for a global existence theory. 
However there are no corresponding uniqueness theorems, and the
framework of measure-valued solutions  is clearly inadequate to distinguish those solutions which are 
physically relevant, 
and has to be supplemented with further structural conditions on the 
solutions. 
Clearly a minimal requirement for any new concept of solution is 
that it should agree with
the classical solution when such exists, and more generally
it is worthwhile to determine properties of classical solutions which 
carry over to the new class of solutions.

We consider the {\em dissipative  measure-valued solutions}
(see definitions \ref{def1}, \ref{dissmvqel}, \ref{defmvpcel}
and  \ref{defmvsol}), 
which form a sub-class of the measure-valued solutions which satisfy an 
averaged and integrated form of the entropy inequality (which allows for
concentration effects in the $L^p\,,\,p<\infty$ setting). We prove 
that, when a classical solution is present, 
the dissipative measure-valued and the classical solution coincide.
The method of proof is based on the idea of relative entropy and the 
format of weak-strong uniqueness 
that was introduced in the context of conservation laws 
in \cite{dafermos79,daf}. The {\it measure-valued-strong uniqueness}
which we prove here handles both
oscillations and concentrations, and it is a further consequence of
the method of proof that when a classical solution exists a dissipative
measure-valued solution does not admit concentrations in the entropy.
To carry out this generalization one needs to account for 
concentrations in the approximating sequence as in
\cite{dipm} and \cite{ab97}. 
For present purposes however we do not need the general
representation of concentrations obtained in these articles, 
because we only consider concentration effects for a single function - the 
entropy which appears in the definition of dissipative solution.
In appendix \ref{concen} we provide a completely elementary derivation
of the {\it Young measure with concentration}  representation of the weak limit
of this function, see \eqref{conc}.

The second issue we study in section \ref{lpb} is the role of  (the measure-valued form of) 
the entropy inequality and the
sense in which {\it entropic measure-valued solutions} assume the initial data. 
Several authors have studied the problem of the initial trace of solutions
for conservation laws, starting with 
\cite{diperna} and then  \cite{chenrascle}, \cite{vasseur}
(using genuine nonlinearity)
and \cite{panov} (exploiting the entropy inequality).
We show that when the Young measure
associated to the family of initial data is a Dirac mass a time average 
of a dissipative measure-valued solution converges strongly 
to the  initial data (see theorem \ref{initialtrace}). 
This result, which extends the
obervations of DiPerna in \cite[section 6(e)]{diperna} to an
$L^p$ context where there is the possibility of
the development of concentrations which has to be eliminated,
represents another noteworthy consequence of the convexity of the entropy.

The  relative entropy method used here to prove  measure-valued-strong uniqueness  provides a clean and quick
proof of strong convergence of approximation schemes to conservation laws 
in the time regime in which  the conservation laws admit classical solutions.
To explain this recall that a conservative view of measure-valued solutions
is that they provide an efficient way of encoding some  properties of
weakly convergent approximating sequences to a system of equations.
Once an approximation scheme is established which is stable, 
in the precise sense 
that it generates a dissipative measure-valued solution,
measure-valued-strong uniqueness automatically implies
strong convergence, {\it without energy concentration}, of the
approximating sequence  to the classical solution.
We illustrate this aspect by considering a lattice approximation 
of the equations of elasticity (in one space dimension) 
by a system of point masses connected by nonlinear springs, 
and prove strong convergence of the spring-mass system to the equations of 
one-dimensional elastodynamics in the continuum limit 
(as long as the latter admits a classical solution).

After the completion of this work we became aware of a recent article 
by Brenier-DeLellis-Sz\`ekelyhidi \cite{bds} in which weak-strong uniqueness 
is proved  for measure-valued solutions of the Euler equations. 
Although the focus of our article is a
different system of equations, with specific intrinsic features - 
notably the lack of uniform convexity and the 
embedding into the enlarged system \eqref{sys1}-\eqref{sys2} via
the null Lagrangians - 
there is overlap both in terms of general ideology and more
specifically of the  material in section \ref{linf} on conservation
laws with $L^\infty$ bounds.
Nevertheless,  we retain this material for explanatory purposes.

The article is organized as follows: in section \ref{relenmv} we introduce
the problem and then in section 
\ref{linf} we perform the basic relative entropy computation 
at the level of a system of conservation laws with $L^\infty$ bounds 
for an approximating sequence, and deduce measure-valued-strong uniqueness.
(theorem \ref{wk}).
Then in section \ref{qel} we generalize to handle the situation that
the approximating sequence is only bounded in $L^2$: we study the
quasi-linear wave equation with convex stored energy satisfying 
quadratic growth conditions above and below, and show how to handle
the possibility of concentrations using the material in appendix
\ref{concen}.
In section \ref{pel} we recall the
global existence of measure-valued solutions for polyconvex elastodynamics
from \cite{dst2} and show that they are dissipative (where the relevant
entropy is the energy, re-interpreted as the convex entropy for the
enlarged system \eqref{sys1}-\eqref{sys2}). We then show that the relative
entropy computation can be performed for this system and prove
measure-valued-strong uniqueness (theorem \ref{pcmvs}). 
In section \ref{lpb} we discuss general systems of
conservation laws with $L^p$ bounds, first extending 
measure-valued-strong uniqueness
to the $L^p$ case in theorem \ref{wkp}
and then proving theorem \ref{initialtrace}
on the strong attainment of the initial data. Finally section \ref{continuum}
is on the lattice-continuum limit for one dimensional elastodynamics.
 
As a final comment, the embedding of polyconvex elastodynamics 
into \eqref{daf1}-\eqref{daf2} 
notwithstanding, 
theorem \ref{pcmvs} is not a consequence of theorem
\ref{wkp} on general systems of conservation laws: both the statement
of the hypotheses for and the proof of theorem \ref{pcmvs} make use of
specific structural features of polyconvexity and the proof requires
the weak  continuity of the null Lagrangians.

\section{Relative entropy for measure-valued solutions}\label{relenmv}

Consider the system of conservation laws,
\begin{equation}
\label{conlaw}
\del_t v + \div_x \,  f(v) = 0 \, ,
\end{equation}
where ${{v}}=({{v}}_1,\dots {{v}}_n)$  are functions of $x=(x_1,\dots x_d)\in\R^d$ and $t\geq 0$.
Attempts to prove an existence theorem for \ref{conlaw} typically involve
the study of a sequence of functions $v^\e$ which are solutions
of an approximating problem
\begin{equation}
\label{appcl}
\del_t v^\eps + \div_x \, f(v^\eps) = \cp_\eps
\end{equation}
where $\cp_\eps \to 0$ in distributions. Uniform bounds for the
sequence are typically a consequence of an entropy inequality 
for the appoximating problem:
\begin{equation}
\label{appentropy}
\del_t \eta (v^\eps)  + \div_x \, q(v^\eps) \le \cq_\eps
\end{equation}
with again $\cq_\eps \to 0$ in distributions. 
Typically \eqref{appentropy} provides
the available uniform bounds, 
$
\sup_{\e,t}\int\eta(v^\e(x,t))\,dx<\infty\,,
$
for the sequence of approximate solutions.
In the limit such an approximation procedure typically 
yields a measure-valued solution verifying
a measure-valued version of the entropy inequality.
One technical difficulty arising here however
is that classical Young measures represent weak limits of functions of growth strictly less than that of $\eta$ but  are insufficient 
to represent the weak limit of $\eta$ itself.
The class of Young measures has to be adapted to reflect the 
representation of the weak limits 
of the entropy function in the presence of concentrations.
We present a self-contained development of a technical  tool designed to address this difficulty in 
appendix \ref{concen}, see \eqref{conc}.
The concentration measure developed there (see \eqref{conc}) will
be incorporated in the definition of the class of dissipative 
measure-valued solutions which are studied in this article. 

In this section we explain in the context of two model problems
how to prove that, in the presence of a classical
solution, a dissipative measure-valued solution with the same initial data
necessarily agrees with that classical solution 
({\it measure-valued-strong uniqueness}).
The presentation is split into two: in section \ref{linf},  
in the presence of  uniform $L^\infty$ bounds, 
classical Young measures are used for the definition of measure-valued 
solution and the basic relative entropy computation
(\cite[Section 5.2]{daf}) is shown to extend
to the measure-valued situation, yielding the
proof of theorem \ref{wk}.
In section \ref{qel}, we take up a model problem for 
the equations of elastodynamics: the quasi-linear wave equation with 
convex quadratic stored energy, where the appropriate stability 
framework involves uniform $L^2$ bounds. 
There, the tool of Young measure with  energy concentration 
developed in appendix \ref{concen} is used to define the appropriate 
notion of {\em dissipative} measure-valued solution,
and this is then used to prove theorem \ref{qlu} on measure-valued-strong
uniqueness in the presence of energy concentration.

\subsection{Conservation laws with $L^\infty$ bounds}
\label{linf}

Consider the system \eqref{conlaw} written in coordinate form, 
\beq
\la{c}
\frac{\partial {{v}}_j}{\partial t}+\frac{\partial {f}_{j\alpha}}{\partial x_\alpha}=0 , , 
\eeq
where latin indices $i,j,k\dots$ are used for the target and greek indices $\alpha,\beta\dots$ for the 
domain. The summation convention will be used throughout. To avoid inessential issues, we will work in the spatially periodic case and spatial integrals will be over the fundamental domain of periodicity $Q=(\R/2\pi\Z)^d$. 
We write $Q_T=Q\times [0,T)$ for $T\in[0,+\infty)$ and $\overline{Q}_T = Q\times [0,T]$.

We assume that \eqref{c} is endowed with an entropy - entropy flux pair $\eta - q$, that is, it is equipped with
an additional conservation law
\beq
\la{ee}
\frac{\partial \eta}{\partial t}+\frac{\partial q_{\alpha}}{\partial x_\alpha}=0 \, ,
\eeq
and that the entropy function $\eta$ is convex.
Then, $\eta - q$ satisfy the consistency equations
\beq
\la{ed}\frac{\partial\eta}{\partial {{v}}_j}
\frac{\partial {f}_{j\alpha}}{\partial {{v}}_i}=\frac{\partial q_\alpha}{\partial {{v}}_i} \, ,
\eeq
or equivalently
\beq
\la{edc}
\frac{\partial^2\eta}{\partial {{v}}_k\partial {{v}}_j}
\frac{\partial {f}_{j\alpha}}{\partial {{v}}_i}
=
\frac{\partial^2\eta}{\partial {{v}}_i\partial {{v}}_j}
\frac{\partial {f}_{j\alpha}}{\partial {{v}}_k}.
\eeq
All functions $f,\eta,q$ are assumed $C^2$ and we assume
positivity of the Hessian matrix $\nabla^2 \eta$ (which implies
strict convexity of $\eta$).

\begin{definition} \la{def1}
\rm
Let $\n = \{ \n_{x,t} \}_{ \{ (x,t) \in \barQT \} }$ be a parametrized 
family of probability measures that are all supported within a compact 
subset $D\subset \R^n$, and with the property that 
for all continuous $f:\R^n\to\R$
$$
\langle\,\n,\,f\,\rangle=\langle\,\n_{x,t},\,f\,\rangle=\int f(\lambda)\,d\n(\lambda)
$$
is a measurable function of $(x,t)$. 
\begin{itemize}
\item[(i)] The pair $({{v}},\n)$ is a  {\em measure-valued solution} of \eqref{c} 
with initial values ${{v}}_0(x)$, if it 
verifies ${{v}}=\int\lambda\, d\n(\lambda)\in L^\infty(dxdt)$ and
\beq
\la{wc}
\iint\biggl[\frac{\partial\psi_i}{\partial t} {{v}}_i+
\frac{\partial\psi_i}{\partial x_\alpha}\langle\n, {f}_{i\alpha}\rangle\biggr]dxdt
+\int\psi_i(x,0){{v}}_{0,i}(x)dx=0,
\eeq
for any test functions $\psi=\psi(x,t) \in C^1_c  ( Q_T )$.

\item[(ii)]
 It will be called an {\it entropic} measure-valued solution of \eqref{c} if,
in addition, for {\em non-negative} test functions, $\psi \in C^1_c  \big ( Q_T \big )$ with
$\psi \geq 0$, there holds:
\beq
\la{wec}
\iint\biggl[\frac{\partial\psi}{\partial t} \langle\n,\eta\rangle+
\frac{\partial\psi}{\partial x_\alpha}\langle\n, q_{\alpha}\rangle\biggr]dxdt
+\int\psi(x,0)\eta({{v}}_{0}(x))dx\geq 0.
\eeq

\item[(iii)]
  It will be called a {\em dissipative} measure-valued solution
if this inequality holds only for non-negative test functions $\psi(x,t)=\theta(t)$ depending solely on time,  i.e. if
\beq
\la{weni}
\iint\frac{d\theta}{d t} \langle\n,\eta\rangle
dxdt
+\int\theta(0)\eta({{v}}_{0}(x))dx\geq 0.
\eeq
for all $\theta \in C^1_c \big ( [0,T) \big )$ satisfying $\theta \geq 0.$
\end{itemize}
\end{definition}

We assume that  there is a classical solution
of \eqref{c} on $\barQT$, to be precise a function  
${{\ov}}\in W^{1,\infty}(\barQT)$
(i.e. a bounded function which is differentiable a.e. 
{\em with bounded derivative}) which
verifies the strong (or classical) versions of \eqref{wc} and \eqref{weni}:
\beq
\la{wcc}
\iint\biggl[\frac{\partial\psi_i}{\partial t}{{\ov}}_i+
\frac{\partial\psi_i}{\partial x_\alpha}{f}_{i\alpha}({\ov})\rangle\biggr]dxdt
+\int\psi_i(x,0){{\ov}}_{0,i}(x)dx=0\,,
\eeq
and
\beq
\la{wenic}
\iint\frac{d\theta}{d t} \eta({\ov})\,
dxdt
+\int\theta(0)\eta({{\ov}}_{0}(x))dx = 0\,,
\eeq
for all test functions $\psi,\theta$ as above.
(Note that \eqref{wenic} is now an equality).
In this circumstance we have the following:
\begin{theorem}
\la{wk}
Let ${{\ov}}\in W^{1,\infty}(\barQT)$ verify \eqref{wcc} and \eqref{wenic}
and let $({v},{\n})$ be a dissipative measure-valued solution verifying
\eqref{wc} and \eqref{weni}.
Assume there exists a compact set $D\subset \R^n$ 
in which ${{\ov}}$ takes its values, 
and assume also that ${{v}}$ takes its
values in $D$, and that $\n$ is supported in $D$. Then there exists
$c_1>0,c_2>0$ such that for $t\in [0,T]$:
\beq\la{cwpp}
\iint\,|\lambda-{{\ov}}|^2\,d\n(\lambda)dx\leq c_1 \left ( \int |{v}_0-{\ov}_0|^2\,dx \right ) \, 
e^{c_2t}\,,
\eeq
and in particular if the initial data agree, ${v}_0={\ov}_0$
then $\n=\delta_{{{\ov}}}$ and ${v}={\ov}$ almost everywhere.
\end{theorem}

\proof
Introduce the {\it relative entropy}
\begin{equation}
\label{defrelen}
\eta_{rel}(\lambda,{\ov}) := \eta(\lambda)-\eta({{\ov}})
-\frac{\partial\eta}{\partial {{v}}_j}({{\ov}})(\lambda_j-{{\ov}}_j) \, ,
\end{equation}
the averaged quantities
\begin{align}
h(\n,{{v}},{{\ov}}) &:= \langle\n,\eta\rangle-\eta({{\ov}})-
\frac{\partial\eta}{\partial {{v}}_j}({{\ov}})({{v}}_j-{{\ov}}_j) \, ,
\\
Z_{k\alpha}(\n,{{v}},{{\ov}}) &:=  \langle\n,{f}_{k\alpha}\rangle-{f}_{k\alpha}({{\ov}})-
\frac{\partial {f}_{k\alpha}}{\partial {{v}}_j}({{\ov}})({{v}}_j-{{\ov}}_j) \, , 
\la{deferr}
\end{align}
and note that, since $\n$ is a probability measure at each $x,t$, it is possible to write
\beq\la{rtm}
h(\n,{{v}},{{\ov}})=\int\Bigl(\eta(\lambda)-\eta({{\ov}})
-\frac{\partial\eta}{\partial {{v}}_j}({{\ov}})(\lambda_j-{{\ov}}_j)\Bigr)\,d\n(\lambda)
=\int\eta_{rel}(\lambda,{\ov}) 
\,d\n(\lambda) \, .
\eeq

Next, using \eqref{c} and \eqref{edc} we calculate that:
\begin{align*}
\frac{\partial}{\partial t}\biggl(
\frac{\partial\eta}{\partial {{v}}_j}({{\ov}})\biggr)=
\frac{\partial{{\ov}}_k}{\partial t}\frac{\partial^2\eta}{\partial {{v}}_k\partial {{v}}_j}
({{\ov}})
&=-\frac{\partial }{\partial x_\alpha}{f}_{k\alpha}({{\ov}})
\frac{\partial^2\eta}{\partial {{v}}_k\partial {{v}}_j}
({{\ov}})\\
&=-
\frac{\partial{{\ov}}_l}{\partial x_\alpha}
\frac{\partial {f}_{k\alpha}}{\partial {{v}}_l}({{\ov}})
\frac{\partial^2\eta}{\partial {{v}}_k\partial {{v}}_j}
({{\ov}})\\
&=-
\frac{\partial{{\ov}}_l}{\partial x_\alpha}
\frac{\partial {f}_{k\alpha}}{\partial {{v}}_j}({{\ov}})
\frac{\partial^2\eta}{\partial {{v}}_k\partial {{v}}_l}
({{\ov}})\,,
\qquad\qquad\hbox{by \eqref{edc}\,.}
\end{align*}
Since this is a bounded function (on account of the 
hypothesis that ${{\ov}}$ is Lipschitz), and referring to
the definition of $Z$ in \eqref{deferr} above, we
deduce that
\beq
\frac{\partial}{\partial t}\biggl(
\frac{\partial\eta}{\partial {{v}}_j}({{\ov}})\biggr)
({{v}}_j-{{\ov}}_j)
+\frac{\partial}{\partial x_\alpha}\biggl(
\frac{\partial\eta}{\partial {{v}}_k}({{\ov}})\biggr)
\Bigl(\langle\n, {f}_{k\alpha}\rangle-{f}_{k\alpha}({{\ov}})\Bigr)
=\frac{\partial{{\ov}}_l}{\partial x_\alpha}
\biggl(\frac{\partial^2\eta}{\partial {{v}}_k\partial {{v}}_l}({{\ov}})\biggr) Z_{k\alpha}.
\la{tt}\eeq

Note that, upon using an approximation argument,
$\psi$ and $\theta$ in \eqref{wc}, \eqref{wec},
\eqref{weni},  \eqref{wcc} and \eqref{wenic} can be taken to be 
Lipschitz functions that vanish for sufficiently large times.
Now choose $\psi(x,t)=\theta(t)\frac{\partial\eta}{\partial {{v}}_j}({{\ov}}(x,t))$
in \eqref{wc} and \eqref{wcc},
subtract them, and then apply \eqref{tt} to get:
$$
\iint\biggl[
\frac{d\theta}{d t} 
\frac{\partial\eta}{\partial {{v}}_j}({{\ov}})({{v}}_j-{{\ov}}_j)
+
\theta
\frac{\partial{{\ov}}_l}{\partial x_\alpha}
\biggl(\frac{\partial^2\eta}{\partial {{v}}_k\partial {{v}}_l}({{\ov}})\biggr) Z_{k\alpha}
\biggr]dxdt
+\int\theta\frac{\partial\eta}{\partial {{v}}_j}({{\ov}})\Bigg|_{t=0}
\bigl[{{v}}_{0,j}(x)-{{\ov}}_{0,j}(x)\bigr]dx=0.
$$
Next, subtract this equation from \eqref{weni}, and also subtract 
\eqref{wenic}, leading to:
\begin{align}
\la{ineq}
\iint\,{\dot{\theta}}\,h\,dxd\tau
&\geq
\iint\theta\frac{\partial{{\ov}}_l}{\partial x_\alpha}
\biggl(\frac{\partial^2\eta}{\partial {{v}}_k\partial {{v}}_l}({{\ov}})\biggr) Z_{k\alpha}
\,dxd\tau\\
&\quad\qquad -\int\theta(0)\Bigl[
\eta({{v}}_0)-\eta({{\ov}}_0)-\frac{\partial \eta}{\partial {{v}}_i}({{\ov}}_0)
({{v}}_0-{{\ov}}_0)_i\Bigr]dx,\nonumber
\end{align}
for non-negative Lipschitz test functions $\theta=\theta(\tau)$.
Now let $\theta(\tau)$ be the non-negative piecewise linear function 
given by
\beq
\theta(\tau)\equiv
\begin{cases}
&1\;\mbox{ when }\;0\leq\tau<t\,,\\
&0\;\mbox{ when }
\tau\geq t+\epsilon\,,\\
&\frac{t-\tau}{\epsilon}+1\quad\mbox{when}\quad
t\leq\tau<t+\epsilon\,.
\end{cases}
\label{bmb}
\eeq
With this choice of $\theta$ \eqref{ineq} reads
\beq\la{utn}
-\frac{1}{\epsilon}\int_{t}^{t+\epsilon} \int h dxd\tau\geq 
\iint\theta(\tau)\frac{\partial{{\ov}}_i}{\partial x_\alpha}
\frac{\partial^2\eta (\bar v)}{\partial {{v}}_k\partial {{v}}_i}Z_{k\alpha}dxd\tau
-\int\Bigl[
\eta({{v}}_0)-\eta({{\ov}}_0)-\frac{\partial \eta}{\partial {{v}}_i}({{\ov}}_0)
({{v}}_0-{{\ov}}_0)_i\Bigr]dx
\eeq
which implies, in the limit $\epsilon\to 0$, that 
\begin{align}
\label{utln}
\int h\,dx  
\leq c\int_0^t\int  \max_{k,\alpha} |Z_{k\alpha}| \,dx\,d\tau
+\int\Bigl[
\eta({{v}}_0)-\eta({{\ov}}_0)-\frac{\partial \eta}{\partial {{v}}_i}({{\ov}}_0)
({{v}}_0-{{\ov}}_0)_i\Bigr]dx 
\end{align}
for $t\in (0,T)$.

Under the working assumption that $\eta$ has strictly positive second
derivative, there exists ${c_0}={c_0}(D)>0$ such that
\beq
\la{lb}
h(\n,{{v}},{{\ov}})\geq{c_0}\int |\lambda-{{\ov}}|^2d\n(\lambda).
\eeq
Notice also that, for some $C = C(D)$,
\beq
\la{until}
\begin{aligned}
|Z_{k\alpha}(\n,{{v}},{{\ov}})|
&=
|\langle\n,  {f}_{k\alpha}(\lambda) -{f}_{k\alpha}({{\ov}})-
\frac{\partial {f}_{k\alpha}}{\partial {{v}}_j}({{\ov}})({{\lambda}}_j-{{\ov}}_j) \rangle | 
\\
&\leq C\int |\lambda-{{\ov}}|^2d\n(\lambda) \, .
\end{aligned}
\eeq
Hence,
\begin{align}
{c_0}\int\, &|\lambda-{{\ov}}|^2\,d\n(\lambda)\,dx 
\leq \int h\,dx  \nonumber
\\
&\leq c\int_0^t\int\int\, |\lambda-{{\ov}}|^2\,d\n(\lambda)\,dx\,d\tau
+c'\int\,|{v}_0-{\ov}_0|^2\,dx\,,
\end{align}
where
$c=c(D,|{{\ov}}|_{W^{1,\infty}})$ and $c'=c'(D)$.
Therefore Gronwall's inequality implies the bound \eqref{cwpp}
and the fact that if ${\ov}_0={v}_0$
then $\int\,|\lambda-{{\ov}}|^2\,d\n(\lambda)\,dx$
is zero at later times, i.e. the measure-valued solution
agrees with the classical solution ${{\ov}}$ almost everywhere. \qed

The above calculation is a measure-valued version of the calculation  
in \cite[Section 5.2]{daf}. In an analogous fashion,  it can 
be carried through for test functions
with more general $x$-dependence to give a measure valued version of 
equation (5.2.6) in that reference, 
but we do not pursue that here.

\subsection{Quasilinear wave equation with convex energy and $L^2$ bounds}
\la{qel} 
In this section we consider the quasi-linear wave equation:
\begin{equation}
\label{mm}
\frac{\partial^2 y}{\partial t^2}=\nabla\cdot S(\nabla y),
\end{equation}
where $y \; : \; {Q} \times {\Real}^+ \to{\Real}^3$ 
and $S$ is the gradient of a strictly convex function
$G:\mat^{3\times 3}\to [0,\infty)$, about which we make the
following hypotheses:
\begin{itemize}
\item[(a1)] $G\in C^3$ and $m|Z|^2\leq D^2G(\hat F)[Z,Z]\leq M|Z|^2$;
\item[(a2)] $G(F)=g_0(F)+\frac{1}{2}|F|^2$ where $\lim_{|F|\to\infty}
\frac{g_0(F)}{1 + |F|^2}=0$.
\item[(a3)]  $\lim_{|F| \to \infty} \frac{ |\nabla_F G(F) |}{1 + |F|^2} = 0$
\item[(a4)] $|D^3 G(F)| \le M$, for some $M > 0$.
\end{itemize}
(We use the summation convention for repeated indices, the norm
$|F|^2=F_{i\alpha}F_{i\alpha}$ and explicitly the second derivative is given by
$D^2G(\hat F)[Z,\tilde Z]=\frac{\partial^2 G(\hat F)}{\p F_{i\alpha}
\p F_{j\beta}}Z_{i\alpha}\tilde Z_{j\beta}$.)
If $y$ is interpreted as a displacement vector this equation could
be regarded as a model for elastodynamics, but the assumption of
convexity is known to be physically unrealistic. We consider a more
realistic model in section \ref{pel}.

A classical solution of \eqref{mm} means a $C^1$ function whose
first derivatives are Lipschitz and verify \eqref{mm}
almost everywhere.
Alternatively,
introducing the notation $v_{i}=\pt y_i$ and $F_{i\alpha}=
\frac{\partial y_i}{\partial x_\alpha}$, 
a classical solution to \eqref{mm} in first order form consists of
a pair $(v,F)$ of Lipschitz functions which solve
\begin{align}
\frac{\partial v_i}{\partial t}&=
\frac{\partial}{\partial x_\alpha}\bigl(
\frac{\partial G}{\p F_{i\alpha}}
\bigr) 
\label{d1}\\
\frac{\partial F_{i\alpha}}{\partial t}&=\frac{\partial v_i}{\partial
x_\alpha}.\label{d2} 
\end{align}
Such a solution
will automatically satisfy the conservation law
\beq\la{cl}
\pt\eta+\partial_\alpha q_\alpha=0\eeq 
where $\eta(v,F) =\frac{1}{2}|v|^2+G(F)$
and $q_\alpha(v,F) =v_i\frac{\p G}{\p F_{i\alpha}}(F)$, and take on the
initial data $v^0(x)=v(0,x)$ and $F^0(x)=F(0,x)$ in the uniform
norm.

\begin{definition} \label{defmvqel} \rm
A {\em measure-valued} solution to \eqref{mm} with initial
data $(v^0,F^0)\in L^2\oplus L^2$ consists of a pair
$(v,F)\in L^\infty(L^2)\oplus L^\infty(L^2)$ and a Young measure
$\n=(\n_{x,t})_{x,t\in \barQT}$ generated by a sequence satisfying \eqref{ube}
such that for $i,\alpha=1,\dots 3$
\begin{align}\label{wd1}
& \int \psi(0,x) v_{i}^0(x) \dx
+ \iint v_{i} \partial_t  \psi \dxdt = \iint  
\left<\n, \frac{\partial G}{\partial F_{i\alpha}}\right> 
\partial_{\alpha} \psi \dxdt
\\
&\int\psi(0,x)F_{i\alpha}^0(x) \dx
+ \iint F_{i\alpha}\partial_t  \psi dxdt=
\iint v_i\partial_\alpha\psi \, dxdt\label{wd2}
\end{align}
for all test functions $\psi=\psi(t,x) \in C^1_c (Q_T) $.
\end{definition}

In order to define a sense in which a measure-valued solution satisfies
the entropy condition \eqref{cl} as an inequality, it is necessary to introduce some method of describing
concentration effects in sequences of approximate solutions. Any natural construction
of a measure-valued solution to \eqref{d1}-\eqref{d2}, e.g. by the
viscosity method or by time-discretization, produces a family of
functions $(v^\epsilon, F^\epsilon)$ of uniformly bounded energy:
\beq\la{ube}
\sup_{\epsilon}\,\sup_{t\geq 0}\,\int\,\eta(v^\epsilon, F^\epsilon)\,dx\,<\,+\infty
\eeq
which are therefore bounded in
$L^\infty(L^2)\oplus L^\infty(L^2)$. Weak limits of such approximate solutions
limit must be represented somehow. For functions of $(v,F)$
of growth at infinity strictly less than quadratic the ordinary Young measure
as developed in \cite{ball88} is sufficient, providing a weakly 
measurable family of probability measures which represent weak limits
of functions of $(v^\epsilon, F^\epsilon)$ which are weakly precompact in $L^1$.
On the other hand, in order to discuss 
the weak limit of quadratic quantities
such as $\eta(v^\epsilon, F^\epsilon)$ 
it is necessary to describe any limiting concentration 
formations in the sequences.   
In appendix \ref{concen} we  introduce a {\em non-negative}
Radon measure $\g$ to measure concentration effects in the energy 
\beq
\g(\psi)=\iint\psi(x,t) \g(dxdt)=\frac{1}{2}\lim_{\epsilon\to 0}
\iint\,\psi\,\bigl(|v^\epsilon|^2-\langle\n_{x,t},|\lambda|^2\rangle
+|F^\epsilon|^2-\langle\n_{x,t},|M|^2\rangle\bigr)\,dxdt\,,
\eeq
for all bounded continuous $\psi$ vanishing for large times,
see \eqref{dg2}. 
(Here ${\n}_{x,t}$ is a probability measure on $\R^3\times\mat^{3\times 3}$, and
we write $(\lambda,M)$ for the coordinates on $\R^3\times\mat^{3\times 3}$
used in the integration with respect to the measure ${\n}$.)
For the class of nonlinear energies $G$ under consideration we will
then have by the Young measure representation (subsequentially):
\begin{align}\label{wd0}
\iint\psi \, \eta (v^\epsilon, F^\epsilon)\,dxdt\,\to\,\iint\psi\bigl(
\langle\n_{x,t},\eta\rangle\,dxdt+\g(dxdt)\bigr)\,,
\end{align}
for all such $\psi$.

The approximate solutions $(v^\epsilon, F^\epsilon)$ are generated 
by families of initial data
\beq\la{wid1}
(v^{\e,0}(x),F^{\e,0}(x))=(v^\epsilon(0,x), F^\epsilon(0,x))\,,
\eeq
converging weakly in $L^2$ to $(v^0(x), F^0(x))$. According to the results of section \ref{el2},
the initial data generate  a Young measure $\m_x$ and an energy concentration measure $\zet(dx)$ 
with the property that (along subsequences)
\beq\la{wid2}
\int\phi(x) g(v^{\e,0}, F^{\e,0})\,dx \to 
\int\phi(x)\langle\m_x, g(\lambda,M)\rangle\,dx
\eeq
for all continuous $\phi$ and subquadratic $g$, and 
\beq\la{wid3}
\int\phi(x) \eta(v^{\e,0}, F^{\e,0})\,dx  \to 
\int\phi(x)\langle\m_x, \eta(\lambda,M)\rangle\,dx + \int \phi(x) \zet(dx)
\eeq
for all continuous $\phi$. 
In this situation we shall refer to
{\em Young measure initial data}  $(v^0,F^0,\m,\zet)$ for brevity.
The important special case that the 
initial data converge strongly 
corresponds to $\zet \equiv 0$ and to the Young measure
$\m_x$ being a Dirac measure.
In the definition of measure valued solutions we think of fixed initial data, or sequences of
data that converge strongly, i.e. $\m_x$ being a Dirac measure. The definition can be easily adjusted to accomodate more general situations.

Assume now that $(v^\e,F^\e)$ is a sequence bounded in 
$L^\infty(L^2)\oplus L^\infty(L^2)$, verifying 
\eqref{wid1}-\eqref{wid3},  which generates the measure-valued solution
verifying \eqref{wd1}-\eqref{wd2}, and the entropy inequality
$$
\int\psi(0,x)\eta(v^{\e,0}, F^{\e,0})dx
+\iint \partial_t\psi\eta(v^\e,F^\e)+\partial_\alpha\psi\,q_\alpha(v^\e,F^\e)
dxdt\geq 0\,,
$$
for $\psi \in C^1(Q_T)$. Taking the limit
$\e\to 0$ and using \eqref{wd0}, \eqref{wid3} (with $\m_x$ a Dirac measure,
$\zet \equiv 0$) motivates 
the following definition of dissipative measure-valued solution:

\begin{definition} \label{dissmvqel} \rm
Given initial data $(v^0,F^0) \in L^2 \oplus L^2$
a {\em dissipative
measure-valued solution with concentration} to \eqref{d1}-\eqref{d2} 
and \eqref{cl}
consists of a pair
$(v,F)\in L^\infty(L^2)\oplus L^\infty(L^2)$, a Young measure
$\n=(\n_{x,t})_{x,t\in \barQT}$ and a non-negative Radon measure 
$\g\in\cM^+(Q_T)$
such that $(v,F,\n)$ is a measure-valued solution verifying 
\eqref{wd1}-\eqref{wd2}, and in addition:
\beq
\la{weceld1}
\iint\,\frac{d \theta}{d t} \, \bigl(\langle\n_{x,t},\eta\rangle
\,
dxdt\,+\g(dxdt)\bigr)
+\int\theta(0) \eta(v^0, F^0) \, dx  \geq 0\,,
\eeq
for all non-negative functions $\theta(t) \in C^1_c ([0,T))$.
\end{definition}
\begin{theorem}
\label{qlu}
Consider a  dissipative measure-valued
solution with concentration to \eqref{d1}-\eqref{d2} 
as just defined, associated to  initial data $(v^0, F^0)$.
\begin{itemize}
\item[(i)]
If $(\hat v, \hat F)\in W^{1,\infty}(\barQT)$ 
is a Lipschitz classical solution with initial 
data $(\hat v^0,\hat F^0)$,
there exist $c_1,c_2>0$ such that for $0\leq t\leq T$:
\beq\la{vbd}
\int\,\langle\n,|\lambda-\hat v|^2
+|M-\hat F|^2\rangle\,dx
\leq 
c_1
\Big ( \int\, |v^0 -\hat v^0 |^2
+|F^0 -\hat F^0 |^2 \,dx   \Big ) \,\, e^{c_2t}\,.
\eeq
\item[(ii)]
If in addition $v^0 = \hat v^0$ and $F^0 = \hat F^0$ almost everywhere,
then $(v,F)=(\hat v,\hat F)$, and 
$\n_{x,t}=\delta_{\hat v(x,t),\hat F(x,t)}$
almost everywhere and the concentration measure $\g$ is null in $Q_T$.
\end{itemize}
\end{theorem}
\proof
Let $(v,F, \n, \g)$ be a dissipative measure-valued solution 
satisfying \eqref{wd1}, \eqref{wd2}
and \eqref{weceld1}. We note that using an approximation argument 
\eqref{wd1}-\eqref{wd2} can be extended to
hold for Lipschitz test functions $\psi$ that vanish for large times:
here we use the assumpion that $\n$ is generated by a sequence verifying
\eqref{ube} which ensures that all quantities in \eqref{wd1}-\eqref{wd2}
lie in $L^1$ under the hypotheses (a1)-(a4) 
and so the bounded convergence theorem applies.
By contrast, \eqref{weceld1} cannot be extended
to this class in the absence of further infomation
about the concentration measure $\g$.

Assume that $(\hat v,\hat F)$ is a classical solution as defined above. 
It will satisfy \eqref{weceld1} as an equality:
\beq
\la{weceldeq}
\iint\,\frac{d \theta}{d t} \, \langle\n_{x,t},\hat\eta\rangle
\,
dxdt\,
+\int\theta(0)\hat\eta_{0}(x)dx =  0\,,
\eeq
where $\hat\eta=\eta(\hat v,\hat F)$ is the energy evaluated along
the solution. Now subtracting from
\eqref{wd1}-\eqref{wd2} the corresponding equations for
the classical solution $(\hat v,\hat F)$, and choosing the
test functions in the resulting equations to be, respectively, 
$\theta(t)\hat v_i$, and 
$\theta(t)\frac{\partial G}{\partial F_{i \alpha}}(\hat F)$, where
$\theta$ is a $C^1$ function of time vanishing for 
sufficiently large times, we obtain the following identity:
\begin{align}
&
\int \theta(0)\,\vh_{i}(0,x)(v_i- \hat v_{i})(0,x) \dx
+\int\theta(0)\,\frac{\partial G}{\partial F_{i\alpha}}
({\hat F}_{i\alpha}(0,x))
\bigl(F_{i\alpha}(0,x)-{\hat F}_{i\alpha}(0,x)\bigr)\,dx\cr
&\qquad+ \iint \Bigl[(v_i-\hat v_{i})\hat v_i
+(F_{i\alpha}-{\hat F}_{i\alpha})\frac{\partial G}{\partial F_{i\alpha}}({\hat F}) 
\Bigr]\,\partial_t  \theta\, dxdt\cr
&\qquad\qquad= 
 \iint\,\theta\,
(\pa \vh_{i} ) \biggl\langle\,\n_{x,t}\,,
\frac{\partial G({M})}{\partial
        {F}_{i\alpha}} - \frac{\partial G({\hat F})}{\partial
        {F}_{i\alpha}}
        - \frac{\partial^{2} G({\hat F})}{\partial
        {F}_{i\alpha}\partial {F}_{j\beta}}({M}_{j\beta} - {\hat F}_{j\beta} ) \biggr\rangle\,
dxdt\; \equiv\mathcal{Q}\,.
\label{dm}
        \end{align}
This calculation is very similar, but simpler,
to one given in full in the next section, and so will
not be written out.

Define the relative entropy as 
\begin{equation}
\label{relentr}
\eta_{rel} (\lambda,{M}; \vh,{\hat F}) \equiv
\frac{1}{2} |\lambda-\vh|^{2} + G({M}) - G({\hat F})
             -\frac{\partial G({\hat F})}{\partial {F}_{i\alpha}}
({M}_{i\alpha} - {\hat F}_{i\alpha})\,,
\end{equation}
and its $t=0$ version as
\begin{equation}
\label{relentr0}
\eta_{rel,0} =\eta_{rel}(\lambda,{M}; \vh^0,{\hat F^0}) \equiv
\frac{1}{2} |\lambda-\vh^0|^{2} + G({M}) - G({\hat F^0})
             -\frac{\partial G({\hat F})}{\partial {F}_{i\alpha}}
({M}_{i\alpha} - {\hat F^0}_{i\alpha})\,.
\end{equation}
Hypotheses (a1) and (a2) guarantee that $\eta_{rel}$ (resp. $\eta_{rel,0}$)
are bounded above and below by multiples of
$|\lambda-\hat v|^2+|M-\hat F|^2$ (resp. $|\lambda-\vh^0|^2+|M-\hat F^0|^2$).
Combining \eqref{weceld1}, \eqref{weceldeq} and \eqref{dm}, we obtain
\beq
\la{weceldr1}
\iint\,{\dot \theta} \,\bigl(\langle\n_{x,\tau},
\eta_{rel}(\lambda,{M}; \vh,{\hat F})\rangle
\,
dxd\tau\,+\g(dxd\tau)\bigr)
+\,\theta(0)\,\int \eta_{rel} (v^0,{F^0}; \vh^0,{\hat F^0}) \,dx  \geq -{\mathcal Q}\,,
\eeq
where $\theta=\theta(\tau) \in C^1_c ([0,T))$.  We would like to choose $\theta$ as in 
\eqref{bmb}, but this is not $C^1$. Therefore we choose a sequence
of  functions $\theta^n \in C^1_c([0,T))$ which are bounded (uniformly in $n$), 
non-increasing and have the property that 
$\dot\theta^n(\tau)\to\dot\theta(\tau)$ for $\tau\neq t,t+\epsilon$.
Since $\dot\theta^n\leq 0$ and $\g\geq 0$, we can discard $\dot\theta^n\g$ in 
the inequality \eqref{weceldr1}.
Referring to \eqref{dm} and substituting in $\theta^n(\tau)$, 
we use assumption (a4). to deduce that there exists 
$C_1=C_1(|\hat v|_{W^{1,\infty}}\,)$ such that for all $n$
\beq
|\mathcal{Q}|\leq C_1\int_0^{t+\epsilon}\int\,
\langle\,\n_{x,t}\,,|M-{\hat F}|^2\rangle\,dx d\tau\,.
\eeq
To take  the limit $n\to\infty$,  note that $\dot\theta^n$ are bounded  and 
so are $\int\,\langle\n_{x,\tau},
\eta_{rel}(\lambda,{M}; \vh,{\hat F})\rangle\, dx$ (by the assumption
on the generation of $\n$ by a sequence verifying \eqref{ube}) 
so that by bounded 
convergence the time integrals converge. We obtain
$$
\frac{1}{\e}\int_t^{t+\e}\int\langle\n_{x,\tau},\eta_{rel}\rangle\,dx\,d\tau
\leq
\int \eta_{rel} (v^0,{F^0}; \vh^0,{\hat F^0}) \,dx
+\,
C_1\int_0^{t+\epsilon}\int\,
\langle\,\n_{x,t}\,,|M-{\hat F}|^2\rangle\,dx d\tau
\,.
$$
Assumptions (a1) and (a2) imply that 
$\langle\n_{x,\tau},\eta_{rel}\rangle\geq \frac{1}{C_2} \langle\n_{x,\tau},|\lambda-\hat v|^2
+|M-\hat F|^2\rangle$ for some $C_2>0$. Consider the function
$\V(\tau)=\int\,\langle\n_{x,\tau},|\lambda-\hat v|^2
+|M-\hat F|^2\rangle\,dx$, which is an averaged variance of the Young 
measure; it satisfies
$$
\frac{1}{\e C_2 }\int_t^{t+\e}\V(\tau)\,d\tau\leq
\int \eta_{rel} (v^0,{F^0}; \vh^0,{\hat F^0}) \,dx
+\,C_1\int_0^{t+\epsilon}\V(\tau)\,d\tau
\,.
$$
Using Lebesgue's theorem, in the limit $\eps \to 0$, $\V(t)$
satisfies
$$\V(t)\leq 
C_2\int \eta_{rel} (v^0,{F^0}; \vh^0,{\hat F^0}) \,dx
\,dx\, +\,C_1 C_2 \int_0^{t}\V(\tau)\,d\tau
\,,
$$
for almost every $t \in (0,T)$. Therefore by Gronwall's inequality
$$
\V(t)\leq 
C_2  e^{C_1 C_2 t}  \int \eta_{rel} (v^0,{F^0}; \vh^0,{\hat F^0}) \,dx \, .
$$
In particular, if the initial data $(v^0, F^0) = (\hat v^0 , \hat F^0)$ a.e.
then the right hand side vanishes, the Young measure has zero variance for almost
every $x,t$, and $\n_{x,t}=\delta_{\hat v(x,t),\hat F(x,t)}$. Going back
to \eqref{weceld1} we deduce that 
$\iint\dot\theta\g(dxdt) \ge 0$ for all $\theta \in C^1_c([0,T))$ with $\theta \ge 0$
and so the concentration measure $\g\geq 0$ 
is in fact identically zero.\qed

\begin{remark}\la{wt} \rm
In writing down \eqref{wd1} in definition \ref{defmvqel} 
the assumption (a3) is used in
order to represent the weak limit of the stress. The situation should be 
contrasted to the Euler
equations, where the flux is of the same order as the energy and the 
description of concentrations 
enters in the definition of measure-valued solutions, 
see Diperna-Majda \cite{dipm}.
\end{remark}

\section{Polyconvex elastodynamics}
\la{pel}

In this section we consider the system of elasticity
\begin{equation}
\label{mainI}
\frac{\partial^2 y}{\partial t^2}=\nabla\cdot S(\nabla y),
\end{equation}
where $y \; : \; {Q} \times {\Real}^+ \to{\Real}^3$ stands for the motion, $F = \nabla y$, $v = \del_t y$,
and $S$ stands for the Piola-Kirchoff stress tensor
obtained as the gradient of a stored energy function, 
$S = \frac{\del W}{\del F}$.  Here we assume that
$W$ is polyconvex, that is  $W(F) = G ( \Phi(F))$ where 
$G:\mat^{3\times 3}\times\mat^{3\times 3}\times \R \to [0,\infty)$
is a strictly convex function and $\Phi(F) = (F ,\cof F, \det F)\in
\mat^{3\times 3}\times\mat^{3\times 3}\times \R$
stands for the vector of null-Lagrangians: $F$, the cofactor matrix $\cof F$
and the determinant $\det F$.

We recall certain formal properties of the equations of polyconvex elasticity 
referring to \cite{qin,daf,dst2} 
for details. Smooth solutions of \eqref{mainI} satisfy the system of 
conservation laws
\begin{align}
\frac{\partial v_i}{\partial t}&=
\frac{\partial}{\partial x^\alpha}\biggl(
\frac{\partial G}{\partial\Xi^A}(\Phi(F))\frac{\partial\Phi^A}{\partial
F_{i\alpha}}(F)\biggr) 
\label{sys1}\\
\frac{\partial\Phi^A(F)}{\partial t}&=\frac{\partial}{\partial
x^\alpha}\biggl(\frac{\partial\Phi^A}{\partial
F_{i\alpha}}(F)v_i\biggr).
\label{sys2} 
\end{align}
In checking this it is necessary to make use of the fact that 
the null-Lagrangians $\Phi(F)$ satisfy
\begin{equation}
\frac{\partial}{\partial x^\alpha}\biggl( 
\frac{\partial\Phi^A}{\partial F_{i\alpha}}(F)\biggr) = 0\,.
\label{nullag}
\end{equation}
Given this, \eqref{sys2} follows from the chain rule and the 
formulae \cite[(2.12-2.13)]{dst2} 
for the derivatives of the null Lagrangians.
In writing the above relations it is
implicitly assumed that $F$ is a gradient 
(which, if it holds initially, is a consequence of 
$\partial_t F = \nabla_x v $, 
 and
this equation is included as the 
first part of \eqref{sys2} since the components of $F$ 
constitute the first nine components of $\Phi(F)$).
Smooth solutions of \eqref{sys1}-\eqref{sys2} automatically satisfy the
conservation of mechanical energy
\begin{equation}
\label{sysenergy}
\del_t \Big ( \frac{1}{2} |v|^2 + G(\Phi(F)) \Big ) 
- \del_\alpha \Big ( v_i\,\frac{\partial G}{\partial\Xi^A}(\Phi(F)) \frac{\partial\Phi^A}{\partial F_{i\alpha}}(F) \Big )
= 0 \, .
\end{equation}

Using these observations
the equations of polyconvex elasticity can be embedded into a
symmetrizable hyperbolic system that determines  the evolution of an
enlarged vector $\Xi=(F,Z,w)$ taking values in
$\mat^{3\times 3}\times \mat^{3\times 3}\times\Real$
and treated as a new dependent variable:
\begin{align}
\frac{\partial v_i}{\partial t}&=
\frac{\partial}{\partial x^\alpha}\biggl(
\frac{\partial G}{\partial\Xi^A}(\Xi)\frac{\partial\Phi^A}{\partial
F_{i\alpha}}(F)\biggr) 
\label{daf1}\\
\frac{\partial\Xi^A}{\partial t}&=\frac{\partial}{\partial
x^\alpha}\biggl(\frac{\partial\Phi^A}{\partial
F_{i\alpha}}(F)v_i\biggr).
\label{daf2} 
\end{align}
Smooth evolutions of this system
preserve the constraints $\Xi^A=\Phi^A(F)$. Moreover, the enlarged system admits the
strictly convex entropy:
\begin{equation}
\eta(v,F,Z,w)=\frac{1}{2}|v|^2+G(F,Z,w)\,,
\label{entropy}
\end{equation}
with corresponding flux 
\beq\la{ef}
q_\alpha=v_i\,\frac{\partial G}{\partial\Xi^A}(\Xi)
\frac{\partial\Phi^A}{\partial
F_{i\alpha}}(F)\,.
\eeq

 We now discuss the various notions of solutions.
A strong (or classical) solution is a $W^{2,\infty}$
function which satisfies \eqref{mainI}; its derivatives automatically
verify \eqref{sys1}-\eqref{sys2} and the
 strong form of the conservation of energy \eqref{sysenergy}.
A weak entropy solution is a weak solution of \eqref{mainI}
which verifies \eqref{sysenergy} as an inequality. In order to make sense 
of the weak forms
the integrability of all quantities which appear has to be guaranteed. 

The notion of measure valued solution that we use is motivated by the 
form of the extended system
\eqref{daf1}-\eqref{daf2} and the existence theory of measure-valued 
solutions developed in \cite{dst2}.
A measure valued solution  will consist of a map $y \; : \; {Q} \times {\R}^+ \to{\Real}^3$, 
with distributional derivatives
$F = \nabla y \in L^\infty ( L^p)$, $v = \del_t y \in L^\infty (L^2)$, and 
a Young measure $\n=(\n_{(x,t)})_{(x,t)\in\barQT}$ 
generated by a sequence satisfying
$$
\sup_{\e, t}\int\,\eta( v^\eps , F^\eps, Z^\eps, w^\eps )\,dx\,<\infty
$$
which represents weak limits in the following way:
\begin{equation}
\label{defym}
\begin{aligned}
\hbox{wk-}\lim_{\epsilon\to 0}
&\;f( v^\eps , F^\eps, Z^\eps, w^\eps )  = \int f(\lav, \laX )  d\n_{(x,t)} (\lav, \laX)
 \\
&\qquad \forall \;\mbox{continuous}  \; f=f(\lav, \laX) \; 
\mbox{with }\; 
\lim_{|\lav| + |\laX| \to \infty} 
\frac{ f(\lav, \laX) }{ \frac{1}{2}|\lav|^2 + G(\laX) } = 0 
\end{aligned}
\end{equation}
where  $\lav \in \R^3$, $\laX = (\laF, \laZ, \law) \in 
\mat^{3\times 3}\times \mat^{3\times 3}\times\Real=\R^{19}$.
The Young measure is connected with the map $y$ through the requirements that
(almost everywhere) 
\begin{equation}
\label{defavg}
F = \langle \n, \laF \rangle  \, , \quad  v = \langle \n , \lav \rangle \, , \quad
\Xi = \langle \n, \laX \rangle \,.
\end{equation}
The action of the Young measure is well defined on all functions that grow slower than the energy
norm. This is the natural framework under the existence of energy norm bounds. With this in mind we define:

\begin{definition}
\label{defmvpcel}
A {\em measure-valued solution} to \eqref{mainI} consists of a map
$y$, with distributional time and space derivatives 
$(v, F) \in L^{\infty}(L^2) \oplus L^{\infty}(L^p) $
and a Young measure $\n=(\n_{x,t})_{x,t\in \barQT}$ as just described, 
such that for $i=1,\dots 3$
\begin{align}
\del_t v_i - \del_\alpha 
\big\langle \n, \frac{\partial G}{\partial\Xi^A}  ( \laX ) \, \frac{\partial\Phi^A}{\partial F_{i\alpha}} (\laF) \big\rangle = 0
\label{defmv1}
\\
\intertext{and  for $A=1,\dots 19$}
\del_t \Phi^A (F)   - 
\del_\alpha  \big ( \frac{\partial\Phi^A}{\partial F_{i\alpha}}(F)  v_i  \big ) = 0
\label{defmv2}
\end{align}
in distributions with
\begin{equation}
\label{defmv3}
\Xi  = \Phi ( \langle \n, \laF \rangle ) = \Phi(F) \, .
\end{equation}

The solution is said to be a {\em dissipative measure-valued solution
with concentration} if it is a measure-valued solution which verifies 
in addition:
\beq
\la{defdissmv4}
\iint\,\frac{d \theta}{d t} \,\Bigl(\langle\n,\eta\rangle
+\g\Bigr)
\,
dxdt
+\int\theta(0)\eta_{0}(x)dx\geq 0\,,
\eeq
for all non-negative functions $\theta=\theta(t) \in C^1_c [0,T)$ with $\theta \ge 0$. Here
$\eta_0$ means the entropy $\eta$ evaluated on the initial data and
$\g$ is the non-negative concentration measure defined in section \ref{elp}.
\end{definition}

The measure-valued solution satisfies 
the momentum equation \eqref{daf1} in the averaged (with respect to the Young measure) sense, 
but the constraint equation  \eqref{daf2} in the {\rm classical weak} sense.  
This is due to the 
weak continuity of the null-Lagrangians 
(see \cite{ball77}, \cite[lemma 3]{dst2}) and the
weak continuity of the transport identities \eqref{sys2} which
follows from the equation $\del_t F = \nabla v$ for 
functions $v\in L^\infty (L^2)$, $F \in L^\infty(L^p)$ 
with $p>4$,  \cite[lemmas 4 and 5]{dst2}.

The existence of a measure-valued solution satisfying \eqref{defmv1}-\eqref{defmv3} is proved in 
\cite[Section 3]{dst2}
under the following hypotheses on the function $G$:
\begin{itemize}
\item[(H1)] $G \in C^3(\mat^{3\times 3} \times \mat^{3\times 3} \times \R ; [0,\infty))$ is a strictly convex 
function satisfying for some $\gamma >0$ the bound $D^2 G \ge \gamma > 0$.
\item[(H2)] $G(F,Z,w) \ge c_1 ( |F|^p +  |Z|^q +  |w|^r + 1)  - c_2$ where $p\in (4, \infty), \ \ q, r \in [2,\infty)$.
\item[(H3)] $G(F,Z,w) \le c (  |F|^p +  |Z|^q + |w|^r +1)$ 
\item[(H4)] 
$| \del_F G|^{\frac{p}{p-1}} + |\del_Z G|^{\frac{p}{p-2}} + |\del_w G|^{\frac{p}{p-3} } \le 
C (  |F|^p +  |Z|^q + |w|^r +1)$
\end{itemize}
The function 
\begin{equation}
\bar G = \alpha |F|^p  + \beta |Z|^q + \gamma|w|^r + |F|^2 + |Z|^2 + w^2
\label{exampleG}
\end{equation}
verifies (H1)-(H3). It
will also verify (H4) under the restrictions $p \ge 2q \ge 4 $, $p \ge 3r \ge 6$.

\begin{theorem} Let $G$ satisfy $(H1)-(H4)$.
Given initial data $(v^0, F^0 )\in L^2 \oplus L^p$, $p \ge 4$,
there exists a dissipative {measure-valued} solution to 
 \eqref{defmv1}-\eqref{defdissmv4} in the sense of definition \ref{defmvpcel}.
\end{theorem}

\proof The existence of a measure-valued solution
 is the main theorem in \cite{dst2}. 
 The fact that this solution
satisfies \eqref{defdissmv4} is proved by using 
the Young measure representation with concentration from section
\ref{elp} to take the limit of equation (3.16)
in \cite{dst2}, using the piecewise constant interpolates
$v^h,\xi^h$ defined in (4.3) in \cite{dst2}, which generate the Young measure
$\n$ in the solution. Using these definitions
equation (3.16) in \cite{dst2} implies that
$$
\int_h^\infty\frac{\theta(t+h)-\theta(t)}{h}\int\,\eta(v^h,\xi^h)\,dx\,dt+\frac{1}{h}\int_0^h\,
\theta(t+h)\,dt\,\int\eta(v^h(x,0),\xi^h(x,0))\,dx\geq 0
$$
for all non-negative functions $\theta(t) \in C^1_c ([0,T))$.
We know that $\frac{\theta(t+h)-\theta(t)}{h}\to\dot\theta(t)$ uniformly as $h\to 0$ , 
But since $\int\eta(v^h,\xi^h)\,dx$ is uniformly bounded this implies that
in this limit we can replace $\frac{\theta(t+h)-\theta(t)}{h}$ by  $\dot\theta(t)$ ,
and then applying \eqref{dg2} we obtain \eqref{defdissmv4}.
\qed

The next objective is to prove the measure-valued-strong uniqueness theorem. In fact the
uniqueness theorem applies to a slightly more general class of nonlinearities: we retain
the hypotheses (H1)-(H3) on $G$, but replace (H4) by the (slightly) weaker hypothesis
\begin{itemize}
\item[(H4)$^\prime$] 
$| \del_F G| + |\del_Z G|^{\frac{p}{p-1}} + |\del_w G|^{\frac{p}{p-2} } \le 
o(1) (  |F|^p +  |Z|^q + |w|^r +1)$ \quad \mbox{ where $o(1) \to 0$ as $|\Xi| \to \infty$}.
\end{itemize}
A hypothesis like $({\rm H}4)'$ is necessary in order to represent the weak limit of the Piola-Kirchhoff stress 
$g_{i\alpha}=\frac{\partial G}{\partial\Xi^A} \big ( \Xi \big ) \, 
\frac{\partial\Phi^A}{\partial F_{i\alpha}}(F)$.
To this end notice that
\begin{align}
\frac{|g_{i \alpha}|}{G(\Xi)} &= \frac{1}{G(\Xi)}  \big |\frac{\partial G}{\partial\Xi^A} \big ( \Xi \big ) \, \frac{\partial\Phi^A}{\partial F_{i\alpha}}(F) \big |
\notag
\\
&\le 
\frac{ | \del_F G | + |\del_Z G| |F| + |\del_w G| |F|^2 }{ |F|^p + |Z|^q + |w|^r + 1}
=
o(1)
\qquad \quad
\mbox{as} \;  |\Xi| \to \infty \, .
\label{pkgrowth}
\end{align}
The last inequality follows from $({\rm H}4)'$ and 
Young's inequality $ab \le \frac{1}{p} a^p + \frac{1}{p'} b^{p'}$, $a,b \ge 0$,
$\frac{1}{p} + \frac{1}{p'} = 1$. By \eqref{pkgrowth} and \eqref{defym} 
the average Piola-Kirchhoff stress $< \n , g_{i \alpha}>$ is then a
well defined locally integrable function which is the weak $L^1$ limit of
$g_{i\alpha}$ evaluated along an approximating sequence. 
As an example notice that the function $\bar G$ in \eqref{exampleG} will satisfy $({\rm H}4)'$ provided
$p > q \ge 2$ and $p > 2r \ge 4$.
We prove:

\begin{theorem} 
\la{pcmvs}
Let $G$ satisfy $(H1)-(H3)$, $(H4)'$ and let  $(y, \n, \g)$ be a dissipative measure-valued 
solution in the sense of definition  \ref{defmvpcel}.
If the initial data equal those of a 
Lipschitz bounded solution $(\hat v,\hat F)\in W^{1,\infty}(\barQT)$: 
$$(v(x,0),\Xi(x,0))=(\hat v(x,0),\Phi(\hat F(x,0)))$$
then $\g$ is zero, $(v,\Xi)=(\hat v,\Phi(\hat F))$ 
and $\n=\delta_{\hat v,\Phi(\hat F)}$.
\end{theorem}

\proof The proof is based on a generalization of the relative entropy computation
to the polyconvex case. Let $(y, \n)$ the measure-valued solution with $v$, $\Xi$ as in \eqref{defavg},
and let $\hat v$, $\hat \Xi := \Phi(\hat F)$ be the Lipschitz solution satisfying \eqref{sys1}-\eqref{sys2}.
As explained in section \ref{qel} we may take the test functions in \eqref{defmv1}
and \eqref{defmv2} to be Lipschitz functions which vanish for large time.
To start with
subtract the weak form of the equations of motion for the measure-valued  
and the Lipschitz
solutions: for $i=1,\dots 3$
\begin{align}\la{w1}
&\int \psi(x,0)(v_i- \hat v_{i})(x,0) \dx
+ \iint (v_i-\hat v_{i}) \partial_t  \psi \dxdt \\
&\qquad=\iint  
\biggl(\left<\n, \frac{\partial G}{\partial\Xi^A}\big( \laX \big) 
\frac{\partial\Phi^A}{\partial F_{i\alpha}} (\laF) \right> 
-\frac{\partial G}{\partial\Xi^A}  \big( \hat\Xi \big) 
\frac{\partial\Phi^A}{\partial F_{i\alpha}}({\hat F})\biggr)
\, \partial_{\alpha} \psi\, \dxdt\notag
\intertext{and  for $A=1,\dots 19$}
\begin{split}
&\int\psi(x,0)\bigl(\Xi^A(x,0)-\hat\Xi^A(x,0)\bigr)\,dx
+ \iint(\Xi^A-\hat\Xi^A) \partial_t  \psi dxdt\\
&\qquad=
\iint \biggl(\frac{\partial\Phi^A}{\partial
F_{i\alpha}}({F})v_i
-\frac{\partial\Phi^A}{\partial
F_{i\alpha}}({\hat F})\hat v_i\biggr) \partial_\alpha\psi \, dxdt
\end{split}\la{w2}
\end{align}
where $\psi$ is a Lipschitz test function that vanishes for 
sufficiently large times. Now choose $\psi$ in \eqref{w1}
to be $\theta(t)\hat v_i$, and in \eqref{w2} to be
$\theta(t)\frac{\partial G}{\partial\Xi^A}(\Phi(\hat F))$, where
$\theta \in C^1_c([0,T))$. Adding the resulting equations 
leads to the identity:

\begin{align}
\begin{split}
&
\int \theta(0) \Big [ \vh_{i}(x,0)(v_i- \hat v_{i})(x,0) 
+  \big ( \frac{\partial G}{\partial\Xi^A}( \hat\Xi^A)
\bigl(\Xi^A-\hat\Xi^A\bigr) \big ) (x,0) \Big ] \, dx
\\
&\quad+ \iint \Bigl[(v_i-\hat v_{i})\hat v_i
+(\Xi^A-\hat\Xi^A)\frac{\partial G}{\partial\Xi^A}(\hat\Xi) 
\Bigr]\,\partial_t  \theta\, dxdt
\\
&\quad= 
 - \iint   \Bigg [ 
 (v_i-\hat v_{i})\pt\hat v_i
+ \bigl(\Xi^A-\hat\Xi^A\bigr)
\pt\bigl(\frac{\partial G}{\partial\xh^A}(\hat\Xi) \bigr)
\notag
 - \partial_{\alpha}\vh_{i} \,  
 \left< \n, 
\frac{\partial G}{\partial\Xi^A}(\laX)
\frac{\partial\Phi^A}{\partial F_{i\alpha}}(\laF) \right> 
\\
&\qquad\qquad+ \partial_{\alpha}\vh_{i}
\frac{\partial G}{\partial\Xi^A}(\hat\Xi)
\frac{\partial\Phi^A}{\partial F_{i\alpha}}({\hat F}) 
-\partial_\alpha
\bigl(\frac{\partial G}{\partial\Xi^A}(\hat\Xi)\bigr)\,
\biggl(\frac{\partial\Phi^A}{\partial
F_{i\alpha}}({F})v_i
-\frac{\partial\Phi^A}{\partial
F_{i\alpha}}({\hat F})\hat v_i\biggr)\, \Bigg ]
\theta \, dxdt
\end{split}
\end{align}

We now calculate, using the fact that $(\hat v,\hat F)$ is a classical
solution of \eqref{daf1}-\eqref{daf2}, and obtain the following identities
for the quantity in square brackets:
\begin{align}
I &:= (\pt\hat v_i) (v_i-\hat v_{i})
+\pt\Big(\frac{\partial G}{\partial\xh^A}(\hat\Xi) \Big)
\bigl(\Xi^A-\hat\Xi^A\bigr)
\nonumber
\\
&\qquad
- \partial_{\alpha}\vh_{i}\,
\Big ( \Big \langle  \n , 
\frac{\partial G}{\partial\Xi^A}(\laX)
\frac{\partial\Phi^A}{\partial F_{i\alpha}}(\laF) \Big\rangle
-\frac{\partial G}{\partial\Xi^A}(\hat\Xi)
\frac{\partial\Phi^A}{\partial F_{i\alpha}}({\hat F}) \Big )
\nonumber
\\
&\qquad 
-\partial_\alpha
\bigl(\frac{\partial G}{\partial\Xi^A}(\hat\Xi)\bigr)\,
\biggl(\frac{\partial\Phi^A}{\partial F_{i\alpha}}({F})v_i
-\frac{\partial\Phi^A}{\partial F_{i\alpha}}({\hat F})\hat v_i\biggr)\, 
\nonumber
\\
 &=  - (\pa \vh_{i} ) 
 \biggl [ \Bigl\langle\,\n\,,  \frac{\partial G}{\partial \Xi^{A}} (\laX) 
    \frac{\partial \Phi^{A}}{\partial F_{i\alpha}} (\laF) \Bigr\rangle 
        -\frac{\partial G}{\partial\Xi^A}(\hat\Xi)
      \frac{\partial \Phi^{A}(\fh)}{\partial  F_{i\alpha}} 
      \nonumber
 \\
        &\qquad \qquad \qquad\qquad\qquad\qquad  - \frac{\partial^{2} G}{\partial
        \Xi^{A}\partial \Xi^{B}} (\hat\Xi))      
        \frac{\partial \Phi^{A}(\fh)}{\partial
        F_{i\alpha}} (\Xi^{B} - \hat\Xi^{B} ) \biggr ]
        \nonumber\\
&\qquad
- \pa  \bigl ( \frac{\partial G }{\partial \Xi^{A}} (\hat\Xi) \bigr ) 
\left ( \frac{\partial \Phi^{A}}{\partial F_{i\alpha}} (F) \, v_{i} 
   - \frac{\partial \Phi^{A}}{\partial F_{i\alpha}} (\fh) \, \vh_{i}
    - \frac{\partial \Phi^{A}}{\partial F_{i\alpha}} (\fh) \, (v_{i} -\vh_{i}) \right )
        \nonumber
\\
        &=  -  (\pa \vh_{i} ) \frac{\partial \Phi^{A}(\fh)}{\partial F_{i\alpha}}
         \biggl\langle\,\n\,,
         \frac{\partial G}{\partial \Xi^{A}} (\laX)
         -\frac{\partial G}{\partial\Xi^A}(\hat\Xi)
         - \frac{\partial^{2} G}{\partial \Xi^{A}\partial \Xi^{B}}(\hat\Xi)    
         (\laXB - \hat\Xi^{B} )
        \biggr\rangle\,
        \nonumber
\\
       &\quad  - (\pa \vh_{i} ) 
       \biggl\langle\,\n\,,
        \Big (  \frac{\partial G}{\partial \Xi^{A}} (\laX)
         -\frac{\partial G}{\partial\Xi^A}(\hat\Xi) \Big )
        \Big ( \frac{\del \Phi^A}{\del F_{i \alpha}} (\laF) - \frac{\del \Phi^A}{\del F_{i \alpha}} (\hat F) \Big )
        \biggr\rangle\,
        \nonumber
\\
&\qquad - \pa  \bigl ( \frac{\partial G}{\partial \Xi^{A}} (\hat\Xi)\bigr )
\Bigl ( \frac{\partial \Phi^{A}(F)}{\partial
        F_{i\alpha}} - \frac{\partial \Phi^{A}(\fh)}{\partial F_{i\alpha}}\Bigr ) (v_i - \hat v_i)
\nonumber
\\
&\qquad\quad
 - (\pa \vh_{i} ) \frac{\partial G}{\partial\Xi^A}(\hat\Xi) 
   \biggl\langle\,\n\,,
 \frac{\del \Phi^A}{\del F_{i \alpha}} (\laF) - \frac{\del \Phi^A}{\del F_{i \alpha}} (\hat F) \biggr\rangle
    \nonumber
\\
&\qquad\qquad - \pa  \bigl ( \frac{\partial G(\hat \Xi) )}{\partial \Xi^{A}}\bigr ) \hat v_i
        \Bigl ( \frac{\partial \Phi^{A}(F)}{\partial
        F_{i\alpha}} - \frac{\partial \Phi^{A}(\fh)}{\partial
        F_{i\alpha}}\Bigr )
\label{lastiden}
\end{align}
Using the fact that
$< \n, \frac{\del \Phi^A}{\del F_{i \alpha}} (\laF) > = \frac{\del \Phi^A}{\del F_{i \alpha}} (F)$ and
the null Lagrangian property \eqref{nullag}, we see that the last two terms can be written as 
a divergence, 
and their contribution integrates to zero. For a test function $\theta \in C^1_c([0,T))$ we obtain:

\begin{align}
&
\int \theta(0) \Big [ \vh_{i}(x,0)(v_i- \hat v_{i})(x,0) + \Big ( \frac{\partial G}{\partial\Xi^A}(\hat\Xi)
\bigl(\Xi^A-\hat\Xi^A\bigr) \Big ) (x,0) \Big ] \,dx
\notag
\\
&\quad+ \iint \Bigl[(v_i-\hat v_{i})\hat v_i
+(\Xi^A-\hat\Xi^A)\frac{\partial G}{\partial\Xi^A}(\hat\Xi) 
\Bigr]\,\partial_t  \theta\, dxdt
= 
 \iint \mathcal{Q} \theta dx dt \, ,
 \label{eqndiff}
\end{align}
where $(-\mathcal{Q})$ stands for the first three terms in \eqref{lastiden},
\begin{equation}
\label{defQ}
\begin{aligned}
\mathcal{Q}    
&=   (\pa \vh_{i} ) \frac{\partial \Phi^{A}(\fh)}{\partial F_{i\alpha}}
         \biggl\langle\,\n\,,
         \frac{\partial G}{\partial \Xi^{A}} (\laX)
         -\frac{\partial G}{\partial\Xi^A}(\hat\Xi)
         - \frac{\partial^{2} G}{\partial \Xi^{A}\partial \Xi^{B}}(\hat\Xi)    
         (\laXB - \hat\Xi^{B} )
        \biggr\rangle\,
\\
       &\qquad  (\pa \vh_{i} ) 
       \biggl\langle\,\n\,,
        \Big (  \frac{\partial G}{\partial \Xi^{A}} (\laX)
         -\frac{\partial G}{\partial\Xi^A}(\hat\Xi) \Big )
        \Big ( \frac{\del \Phi^A}{\del F_{i \alpha}} (\laF) - \frac{\del \Phi^A}{\del F_{i \alpha}} (\hat F) \Big )
        \biggr\rangle\,
\\
&\qquad\quad \pa  \bigl ( \frac{\partial G}{\partial \Xi^{A}} (\hat\Xi)\bigr )
\Bigl ( \frac{\partial \Phi^{A}(F)}{\partial
        F_{i\alpha}} - \frac{\partial \Phi^{A}(\fh)}{\partial F_{i\alpha}}\Bigr ) (v_i - \hat v_i)
\\
&=: Q_1 + Q_2 + Q_3
\end{aligned}
\end{equation}

Defining the relative entropy as 
\begin{equation}
\label{relentrext}
\eta_{rel} (v, \Xi ; \hat v , \hat\Xi ) :=
\frac{1}{2} |v-\vh|^{2} + G(\Xi) - G(\hat\Xi)
             -\frac{\partial G}{\partial \Xi^{A}}(\hat\Xi) \, (\Xi^A - \hat\Xi^A)
\end{equation}
we prove that ${\mathcal Q}$ can be bounded by the averaged relative entropy:
\begin{lemma}
\label{lembds}
Under Hypothesis $(H1)-(H3)$, $(H4)'$, there exists $C=C(|(\hat v,\hat\Xi)|_{W^{1,^\infty}})$ such that
$$
\begin{aligned}
|\mathcal{Q}| &\le
C \langle\n,\eta_{rel}\rangle\,,
\quad
\langle\n,\eta_{rel}\rangle=
\int \eta_{rel}(\lambda_v, \lambda_\Xi ; \hat v , \hat\Xi )\,
\n(d\lambda_v,d\lambda_\Xi)\,.
\end{aligned}
$$
\end{lemma}

{\em Proof of the lemma}.
We start by estimating the term $Q_2$ in  \eqref{defQ}.  
Let $K \subset \R^{19}$ be a compact set containing the values of $\hat \Xi(x,t)$ for $(x,t) \in Q_T$.
We will show that there is a constant $C$ such that for all
$ \laX \in \R^{19}$ and $\hat\Xi \in K \,$ there holds
\begin{equation}
\label{bdI2}
|\cq_2| = \left | \Big (  \frac{\partial G}{\partial \Xi^{A}} (\laX)
         -\frac{\partial G}{\partial\Xi^A}(\hat\Xi) \Big )
        \Big ( \frac{\del \Phi^A}{\del F_{i \alpha}} (\laF) - \frac{\del \Phi^A}{\del F_{i \alpha}} (\hat F) \Big ) \right |
\le C G_{rel} ( \laX ; \hat \Xi) ,
\end{equation}
where
\begin{equation}
\label{bdGrel}
G_{rel} ( \laX ; \hat \Xi) = G(\laX) - G(\hat\Xi) - D_\Xi G(\hat\Xi) \cdot (\laX -\hat\Xi)
\end{equation}

Note that the assumptions $(H1)-(H2)$ imply the lower bound
\beq
\la{lbr}
G_{rel}( \laX ;  \hat\Xi )\geq
\max\{
\gamma (|\laX-\hat\Xi|^2 ,
\alpha (|\laF|^p+|\laZ|^q+|\law|^r+1)- A
\}
\eeq
for some constants $\alpha$, $\gamma$ and $A$  which depend upon $|\Xi|_{L^\infty}$ 
and the constants $c_1, c_2$ appearing in  $(H1)-(H2)$. 

Define now 
${\cal L}_R =\{|\laF|^p+|\laZ|^q+|\law|^r+1 \ge   R \}$ 
with $R$ chosen sufficiently large so that $K \subset ({\cal L}_R)^c$ 
and also 
$$
\alpha (|\laF|^p+|\laZ|^q+|\law|^r+1)- A \geq \frac{\alpha}{2}(|\laF|^p+|\laZ|^q+|\law|^r+1)  
\quad \mbox{ on ${\cal L}_R$}.
$$
For $\laX \in \mathcal{L}_R$ and $\hat\Xi \in K$ we have upon using Young's inequality,
hypothesis $(H4)'$, selecting $R$ sufficiently large, and using \eqref{lbr}  that
\begin{align*}
|\cq_2| &\leq
C\Big [
(1+|\partial_F G(\lambda_\Xi)|)
+(1+|\lambda_F|)(1+|\partial_Z G(\lambda_\Xi)|) +(1+|\lambda_F|^2)(1+|\partial_w G(\lambda_\Xi)|)
\Big]
\\
&\le  \frac{\alpha}{4}  |\laF|^p + C_\alpha \Big ( 
| \del_F G| + |\del_Z G|^{\frac{p}{p-1}} + |\del_w G|^{\frac{p}{p-2} } \Big )
\\
&\le
\frac{\alpha}{2}\big( |\laF|^p+|\laZ|^q+|\law|^r+1 \big)
\\
&\le 
C G_{rel} (\laX ; \hat\Xi) \qquad \laX \in \mathcal{L}_R \, , \; \hat\Xi \in K \, .
\end{align*}
With $R$ now fixed, observe that for $\laX \in (\mathcal{L}_R)^c$
\begin{align*}
|\cq_2| &\leq C_R | \laX - \hat\Xi|^2 
\\
&\le 
\frac{C_R}{\gamma} G_{rel} (\laX ; \hat\Xi) \qquad \laX \in (\mathcal{L}_R)^c \, , \; \hat\Xi \in K \, .
\end{align*}
Therefore, \eqref{bdI2} follows.

The term $Q_1$ is estimated using the bound
\begin{equation}
\label{bdI1}
 |\cq_1| = \big | \frac{\partial G}{\partial \Xi^{A}} (\laX)
         -\frac{\partial G}{\partial\Xi^A}(\hat\Xi)
         - \frac{\partial^{2} G}{\partial \Xi^{A}\partial \Xi^{B}}(\hat\Xi)    
         (\laXB - \hat\Xi^{B} ) \big |
 \le C G_{rel} ( \laX ; \hat \Xi) \quad  \laX \in \R^{19}, \hat\Xi \in K \, ,
\end{equation}
which follows from an argument similar to the derivation of \eqref{bdI2} above (using the
fact from ${\rm H}4'$ that the derivatives of $G$ grow more slowly than $G$ itself at infinity).

Finally the term $Q_3$ is estimated using
\begin{equation}
\label{bdI3a}
|v-\hat v|^2=|\int(\lambda_v-\hat v)d\n|^2\leq
\int|\lambda_v-\hat v|^2d\n\leq C
\int \eta_{rel}(\lambda_v, \lambda_\Xi ; \hat v , \hat\Xi )\,
d\n( \lav,\laX)\,,
\end{equation}
the weak continuity property 
$< \n, \frac{\del \Phi^A}{\del F_{i \alpha}} (\laF) > 
= \frac{\del \Phi^A}{\del F_{i \alpha}} (F)$ 
and the estimation (in the spirit of \eqref{bdI2})
$$
\Big |
\frac{\del \Phi^A}{\del F_{i \alpha}} (\laF) 
- \frac{\del \Phi^A}{\del F_{i \alpha}} (\hat F)
\Big |^2
\le C G_{rel} ( \laX ; \hat \Xi) \quad  \laX \in \R^{19}, \hat\Xi \in K \, .
$$
Combining these we obtain
\begin{align*}
\Big | \frac{\del \Phi^A}{\del F_{i \alpha}} (F) - \frac{\del \Phi^A}{\del F_{i \alpha}} (\hat F) \Big |^2
&=\Big |\int
\big ( \frac{\del \Phi^A}{\del F_{i \alpha}} (\laF) 
- \frac{\del \Phi^A}{\del F_{i \alpha}} (\hat F) \big )
d\n\Big |^2
\notag
\\
&\leq \int
\Big |\frac{\del \Phi^A}{\del F_{i \alpha}} (\laF) 
- \frac{\del \Phi^A}{\del F_{i \alpha}} (\hat F) \Big |^2
d\n
\notag
\\
&\leq C
\int G_{rel}(\lambda_\Xi ; \hat\Xi )\,
d\n(\laX)\,,
\end{align*}
and hence, by \eqref{bdI3a} and Cauchy-Schwarz, 
\beq|{\mathcal Q}_3|\leq C
\int \eta_{rel}(\lambda_v, \lambda_\Xi ; \hat v , \hat\Xi )\,
d\n( \lav,\laX)\,.
\label{bdI3b}
\eeq
The proof of the lemma is completed by refering to 
\eqref{bdI1}, \eqref{bdI2} and \eqref{bdI3b}.
\qed

To conclude, from the definition of the dissipative measure valued solution \eqref{defdissmv4} and the
equations \eqref{eqndiff}, \eqref{defQ}, and lemma \ref{lembds}  we derive the
equation for the relative entropy
\beq
\la{weceldr}
\begin{aligned}
&\iint\,\frac{d \theta}{d t} \, \Big ( \langle\n,\eta_{rel}\rangle \, dxdt + \g (dxdt) \Big )
\\
&\qquad
+\,\theta(0)\,\int\Bigl[\eta_{0}-\hat\eta_0
-\vh_{i}(v_i- \hat v_{i})
-
\frac{\partial G}{\partial\Xi^A}(\xh^A)
\bigl(\Xi^A-\xh^A\bigr)
\Bigr]_{t=0}\,dx\geq -C\int\, \langle\n,\eta_{rel}\rangle \, dxdt\,,
\end{aligned}
\eeq
for all $\theta=\theta(t) \in C^1_c ([0,T))$, $\theta \ge 0$.
The proof can now be completed as in the proof of theorem \ref{qlu}, leading to
the bound
$$
\int\,\langle\n,\eta_{rel}\rangle\,dx\,\bigr|_{t}
\leq 
C_2  e^{C_1 C_2 t}
\int\Bigl[\eta_{0}-\hat\eta_0
-\vh_{i}(v_i- \hat v_{i})
-
\frac{\partial G}{\partial\Xi^A}(\xh^A)
\bigl(\Xi^A-\xh^A\bigr)
\Bigr]_{t=0}\,dx
\, .
$$
This implies the uniqueness assertion in the theorem statement for
appropriate initial data.
\qed

\section{Conservation laws with $L^p$ bounds}
\la{lpb}

In this section we consider a measure-valued solution for 
the system of $n$ conservation laws \eqref{c} in the presence
of $L^p$ bounds for $1< p<\infty$. We first show how to generalize theorem
\ref{wk} on recovery of classical solutions to this case.  We also
discuss the problem of the initial trace, i.e. the sense in which 
a {measure-valued} solution assumes  the initial data. 
In this latter regard we show that the presence of a 
convex entropy yields {\em strong} convergence of the averages
$\frac{1}{\tau} \int_0^\tau v(\cdot , t ) dt $ to the
initial data, thus extending a result of DiPerna \cite{diperna} to the $L^p$ framework.

We assume that \eqref{conlaw} is equipped with an entropy-entropy flux 
pair $\eta - q$ as in section \ref{linf}
with the entropy $\eta$ satisfying  the hypotheses:
\begin{align}
\la{gr0}
\eta \; \mbox{positive and strictly convex}, \; D^2 \eta \ge \gamma > 0
\\
\la{gr3}
\alpha \big (  |\lambda|^p + 1 \big )  - A  \le \eta (\lambda ) \le C \big ( |\lambda|^p + 1 \big )
  \quad \text{$\lambda \in \R^n $}
\end{align}
for some $\alpha, \gamma, A, C > 0$ and for $p \in [2, \infty)$, 
while the flux $f$ in \eqref{c} verifies the growth restriction
\begin{equation}\la{gr1}
\frac{|f(\lambda) | }{\eta (\lambda)} = o(1)  \quad \text{as $|\lambda| \to \infty$}\, .
\end{equation}
The entropy identity provides stability in an $L^p$-framework, $p<\infty$.
In contrast to the $L^\infty$ case treated in section \ref{linf} 
such a  framework permits the development of concentrations in approximating 
sequences, which we describe using the measure ${\g}$ defined in appendix
\ref{concen}.
Using the Young-measure associated to the family $\{ {v}^\eps \}$  and the concentration measure 
${\g} \ge 0$ we have
\begin{align}
g ({v}^\eps) &\rightharpoonup  \langle {\n}_{x,t} , g (\lambda) \rangle \, , \quad \forall  \; 
\text{$g$ continuous s.t. $\lim_{|\lambda| \to \infty} \frac{g(\lambda)}{\eta(\lambda)} = 0$, }
\label{wkg}
\\
\eta({v}^\eps)\, dxdt &\rightharpoonup \langle {\n}_{x,t} , \eta \rangle\, dxdt + {\g}(dx dt)
\label{wketa}
\end{align}
where ${\n}$ and ${\g}$ in \eqref{wketa} are as introduced in 
appendix \ref{concen}.

For the initial data $\{ {v}_0^\eps \}$ of the approximating problem \eqref{appcl} we assume weak convergence to ${v}_0$ in $L^p$ with associated 
Young measure ${\m}_x$, and also allow the development of concentrations in 
$\eta$ described by a  
concentration measure ${\zet}(dx) \ge 0$ such that
\begin{align}
g ({v}_0^\eps) &\rightharpoonup  \langle {\m}_{x} , g (\lambda) \rangle \, , \quad \forall  \; 
\text{$g$ continuous s.t. $\lim_{|\lambda| \to \infty} \frac{g(\lambda)}{\eta(\lambda)} = 0$, }
\label{wkg0}
\\
\eta({v}_0^\eps)\, dx &\rightharpoonup \langle {\m}_{x} , g (\lambda) 
\rangle\, dx + {\zet} (dx) \, .
\label{wketa0}
\end{align}

\begin{definition}
\label{defmvsol}
{\rm A {\em dissipative measure-valued solution with concentration} to 
\eqref{conlaw} 
consists of ${v} \in L^\infty(L^p)$, a Young measure
$(\n_{x,t})_{x,t\in \barQT}$ and a non-negative Radon measure 
$\g\in\cM^+(Q_T)$ such that  
\begin{equation}
\label{eqmv}
\iint  \langle {\n}_{x,t} , \lambda \rangle \psi_t \, dx dt 
+
\iint  \langle {\n}_{x,t} , f(\lambda) \rangle \psi_x \, dx dt 
+
\int {v}_0 (x) \psi (x,0)  \, dx  = 0
\end{equation}
for any $\psi \in C^1_c (Q \times [0,T))$, and
\beq
\la{eqdiss}
\iint  \dot\theta \bigl( \langle{\n}_{x,t},\eta (\lambda) \rangle \,dxdt\,+{\g} (dxdt)\bigr)
+\int\theta(0) \bigl( \langle {\m}_x,\eta\rangle\, dx + {\zet}(dx) \bigr)  \geq 0\,,
\eeq
for all $\theta = \theta(t)  \in C^1_c ([0,T))$ with $\theta \ge 0$.
}\end{definition}

\subsection{Recovery of classical solutions from measure-valued solutions}

We first state the generalization of theorem \ref{wk} in the $L^p$ framework:

\begin{theorem}
\la{wkp}
Let $({v},{\n}, \g)$ be a dissipative measure-valued solution
as in definition \ref{defmvsol}, and suppose that there exists a strong solution 
${{\ov}}\in W^{1,\infty}(\barQT)$ verifying \eqref{wcc} and \eqref{wenic}.
If for the initial data ${\zet}=0$ and ${\m}_x=\delta_{{\ov}_0(x)}$ then 
$\n=\delta_{{{\ov}}}$ and ${v}={\ov}$ almost everywhere on $Q_T$.
\end{theorem}

\proof
The initial calculations are identical to the $L^\infty$ case
in the proof of theorem \ref{wk} up to \eqref{utln}.
Since the support of ${\n}$ is no longer bounded it is necessary to
replace \eqref{until}. This is done as follows: 
define $\eta_{rel}(\lambda,{\ov})$
by \eqref{defrelen} and let $K \subset \R^n$ be a compact set containing the values of
$\bar v(x,t)$ for $(x,t) \in Q_T$. Using \eqref{gr0}, \eqref{gr3}, \eqref{gr1} and an argument
as in the proof of \eqref{bdI2} (see lemma \ref{lembds}), there exists a constant 
$C_1>0$ such that
\beq
\la{high}
\big | f_{k\alpha}(\lambda)-f_{k\alpha}({\ov})-\frac{\partial f_{k\alpha}}{\partial {v}_j}
({\ov})(\lambda_j-{\ov}_j) \big |
\le
C_1  \eta_{rel} ( \lambda ; \bar v) 
\qquad \lambda \in \R^n \, , \; \bar v \in K
\eeq
and hence integrating over $\lambda$ we obtain that
\beq\la{noon}
|Z_{k\alpha}(\n,{{v}},{{\ov}})|
\leq C_1 h(\n,{{v}},{{\ov}})\,,
\eeq
where we use the definitions \eqref{defrelen}-\eqref{deferr}.
This inequality serves as a suitable replacement of \eqref{until}
to complete the transposition of the proof of theorem \ref{wk}
to the $L^p$ setting: under the assumption ${\zet}=0$ there holds
\beq\la{cwp}
\int\, h(\n,{{v}},{{\ov}})\,dx\leq c_1\,\int \eta_{rel}(\lambda,{\ov}_0)d\m(\lambda)
\,dx\,\, 
e^{c_2t}\,,
\eeq
and in particular if ${v}(x,0)={\ov}_0(x)$ and $\m_x=\delta_{{\ov}_0(x)}$ then
${\n}_{x,t}=\delta_{{\ov}(x,t)}$ and ${v}(x,t)={\ov}(x,t)$ for $t> 0$, and
$\g=0$.
\qed

\subsection{On the initial trace of {measure-valued} solutions}\label{secintr}

DiPerna \cite[section 6(e)]{diperna} gave an argument indicating that the
measure-valued version of the entropy condition, 
used in the case of
strict convexity of the entropy,
leads to a strong initial trace for a measure-valued solution 
in the $L^\infty$ setting. 
Below this result is extended to the $L^p$ functional setting, $p<\infty$. 

\begin{theorem}
\label{initialtrace}
Let ${v}$, ${\n}_{x,t}$ and ${\g}(dxdt)$ be a dissipative {measure-valued}
solution with concentration to \eqref{conlaw}. If the Young measure associated with the
data satisfies ${\zet} \equiv 0$ and ${\m}_x = \delta_{{v}_0(x)}$, then as $\tau \to 0 + $
\begin{equation}
\label{strtrace}
\frac{1}{\tau} \int_0^\tau {v}(\cdot , t ) dt \to {v}_0  \, , \quad 
\text{ in $L^p ({Q})$}.
\end{equation}
\end{theorem}

\proof
We first show that as a consequence of the definition of a measure-valued solution 
\begin{equation}
\label{wktrace}
\frac{1}{\tau} \int_0^\tau {v}(\cdot , t ) dt \rightharpoonup v_0  \, , \quad 
\text{weakly in $L^p ({Q})$}\, .
\end{equation}
To achieve this apply \eqref{eqmv} to the test function 
$\psi (x,t) = \varphi(x) \theta (t)$, where
$\varphi \in C^1({Q})$ and 
\beq
\theta(t)\equiv
\begin{cases}
& 1 - \frac{t}{\delta} \quad\mbox{when}\quad
0 \le t\leq \delta \, ,
\\
&0\;\mbox{ when }
\delta \le  t \, .
\end{cases}
\label{testtime}
\eeq
Then we obtain
$$
\begin{aligned}
-\frac{1}{\delta} \int_0^\delta \int_{Q} {v}(x,t) \varphi (x) dx dt 
&+ \int_0^\delta \int_{Q} \langle  {\n}_{x,t} , f(\lambda ) \rangle \varphi (x) \theta (t) dx dt 
\\
&+ \int_{Q} v_0(x) \varphi(x) dx = 0\ .
\end{aligned}
$$
Passing to the limit $\delta \to 0$, we conclude
\begin{equation}
\label{wktracedist}
\lim_{\delta \to 0} \int_{Q}  \left( \frac{1}{\delta} \int_0^\delta  {v}(x,t) dt \right) \varphi (x) dx \to
\int {v}_0(x) \varphi (x) dx\  .
\end{equation}
Since
\begin{equation}
\label{lpbound}
\int_{Q} \left | \frac{1}{\delta} \int_0^\delta  {v}(x,t) dt \right |^p dx 
\le
\frac{1}{\delta} \int_{Q} \int_0^\delta |{v}|^p dx dt
\le 
\| {v} \|_{L^\infty ( L^p ) }
\end{equation}
equation  \eqref{wktracedist}, together with an approximation argument, 
implies that the sequence $\Big \{ \frac{1}{\delta} \int_0^\delta {v}(\cdot , t ) dt \Big \}$
converges weakly to ${v}_0$ in $L^p ({Q})$.

Consider now the functional $I : L^p ({Q}) \to \R$ defined by
$$
I[v] = \int_{Q} \eta (v) dx \, .
$$
Due to the convexity of $\eta$ the functional $I$ is weakly lower semicontinuous. Hence 
\eqref{wktrace} implies
\begin{equation}
\label{wklsc}
\int_{Q} \eta ( {v}_0 (x) ) dx \le \liminf_{\delta \to 0}
\int_{Q} \eta \left ( \frac{1}{\delta} \int_0^\delta {v}(x,t) dt \right ) dx
\end{equation}

Fix  $\theta$ as in \eqref{testtime} and consider a sequence of $C^1$ 
functions $\theta_n\to\theta$ that are monotone decreasing, 
vanish for large $t$, and
satisfy $\theta_n (0) = 1$ and $\dot\theta_n(t) \to \dot\theta(t)$ for
$t\neq 0,\delta$. We apply 
\eqref{eqdiss} to the test functions
$\theta_n$ and use the hypotheses for the initial measure and the property 
${\g} \ge 0$ to obtain
$$
\int_{Q} \eta ({v}_0 (x) ) dx \ge - \iint \frac{d \theta_n}{dt} \langle {\n}_{x,t} , \eta (\lambda) \rangle dx dt \, .
$$
Passing to the limit $n \to \infty$ and then $\delta \to 0$ and using 
$v(x,t)=\int \lambda d{\n}_{x,t} (\lambda)$ and 
Jensen's inequality we conclude
that
\begin{align}
\int_{Q} \eta ( {v}_0 (x) ) dx 
&\ge \limsup_{\delta \to 0}
\frac{1}{\delta} \int_0^\delta \int_{Q} \int \eta ( \lambda)  d{\n}_{x,t} (\lambda)  dx dt \, = \,
\limsup_{\delta \to 0} \int_{Q}\frac{1}{\delta} \int_0^\delta \int \eta ( \lambda)  d{\n}_{x,t} (\lambda)  dt dx
\notag\\
&\ge \limsup_{\delta \to 0}
\int_{Q}   \frac{1}{\delta} \int_0^\delta 
\eta \left (\int \lambda d{\n}_{x,t} (\lambda) \right ) dtdx\,
=\,
\limsup_{\delta \to 0}
\int_{Q}   \frac{1}{\delta} \int_0^\delta 
\eta (v(x,t)) dtdx\,
\notag\\
&\geq\,
\limsup_{\delta \to 0}
\int_{Q} \eta \left (  \frac{1}{\delta} \int_0^\delta 
v(x,t) dt\right ) dx\,.
\la{otherside}
\end{align}
In summary, for the family $\big \{ v^\delta = \frac{1}{\delta} \int_0^\delta {v}(\cdot , t ) dt \big \}$, we have
 $v^\delta \rightharpoonup {v}_0$ weakly in $L^p ({Q})$ and
\begin{equation}
\label{wkcont}
\lim_{\delta \to 0} \int_{Q} \eta \big (  v^\delta (x) \big ) dx 
= 
\int_{Q} \eta ( {v}_0 (x) ) dx \ .
\end{equation}
We claim this implies
\begin{equation}
v^\delta = \frac{1}{\delta} \int_0^\delta {v}(\cdot , t ) dt \to {v}_0  \, , \quad 
\text{ in $L^p ({Q})$}.
\end{equation}
Indeed, by \eqref{lpbound}, the sequence $\{ v^\delta \}$ is uniformly bounded in $L^p({Q})$. The
results of section \eqref{conc} imply that there exists an associated Young measure ${\kap}_x$
and a concentration measure $\ebb(dx) \ge 0 $ such that
\begin{equation}
\label{wkrep}
\eta( v^\delta ) \rightharpoonup \int \eta (\lambda) d{\kap}_x (\lambda) 
+\ebb (dx)
\end{equation}
Now \eqref{wktrace} implies that
$ \int \lambda d{\kap}_x (\lambda)={v}_0(x)$, 
so that by \eqref{wkcont} and \eqref{wkrep} we get
$$
\int_{Q} \int \eta (\lambda) d{\kap}_x (\lambda) \, dx + \int_{Q} \ebb (dx) 
= \int_{Q} \eta ( {v}_0 (x) ) dx
= \int_{Q} \eta \left ( \int \lambda d{\kap}_x (\lambda) \right ) dx\,.
$$
Using Jensen's inequality
\begin{equation}
\label{jens}
\eta \left ( \int \lambda d{\kap}_x (\lambda) \right ) 
\le \int \eta (\lambda) d{\kap}_x (\lambda)
\end{equation}
we conclude that the concentration measure ${\ebb} \equiv 0$, 
and that necessarily \eqref{jens} 
holds as equality. The latter
implies that ${\kap}_x = \delta_{{v}_0(x)}$ and completes the proof of 
\eqref{strtrace}. \qed

\section{Application: one dimensional elastodynamics as the 
continuum limit of a lattice model}
\la{continuum}

Here we investigate a spatially discrete lattice approximation 
to one dimensional
elastodynamics. Apart from interest in the continuum limit, the purpose
is to show that the use of the relative entropy method 
provides an efficient way of proving strong 
convergence theorems for approximation schemes: it is only necessary to 
verify that the approximation scheme generates a dissipative measure-valued
solution.
For simplicity as above  we consider the periodic case 
so that the spatial domain is $Q= \R/ 2\pi\Z$ 
on which are located
$N$ atoms at the points $\{{\poseul}_i(t)\}_{i=0}^{N-1}$, at time $t$,
and continued periodically ${\poseul}_{N+i}(t) = {\poseul}_i(t)+2\pi \  \forall i$ when
convenient.  We assume the existence of an 
equilibrium configuration in which 
the atoms form a one dimensional array
(lattice) in which the $i^{th}$ atom has reference location 
${\poslag}_i =\frac{2\pi i}{N}= \e i $
so they are all separated by a distance 
${\e}\equiv\frac{2\pi}{N}$ from their nearest neighbours on either side.
We write $I^i_\e = \{{\poslag}: {\poslag}_i\leq {\poslag}< {\poslag}_{i+1}\}$ for the intervals into which the 
domain is sub-divided by the reference locations ${\poslag}_i$.

We will 
assume the dynamics is determined by a natural Lagrangian system of
the following form: 
\begin{itemize}
\item
each atom has identical mass 
$\e\rho = \frac{2\pi}{N}\rho$ (so that the  total mass is
$2\pi\rho$), and
the kinetic energy is $T = \frac{1}{2}\e\rho \Sigma_i
\dot{{\poseul}_i}^2$;  
\item
the potential energy is given
by $V= \sum_{i=0}^{N-1} W(\frac{{\poseul}_{i+1} - {\poseul}_i}{\e})$, where
$W$ is a strictly convex $C^3$ function such that
$W''(u)\geq c_0>0$ and 
$
W(u)\geq max(0, c_1|u|^p - c_2)
$
for $c_i > 0 $ , $p\geq 2$ and $u\in\R$;
\item
$\lim_{|u|\to +\infty}\frac{W'(u)}{|u|^p}=0$\;
\item
finally, the Lagrangian
$$
L = T-V = \sum_{i=0}^{N-1} \frac{\e \rho}{2}\dot{{\poseul}^2_i}
-\e W(\frac{{\poseul}_{i+1} - {\poseul}_i}{\e}) .
$$
\end{itemize}
Thus we have the following {\it equation of motion}
\be
\la{eom}
\frac{d}{dt} (\e \rho \dot{{\poseul}}_i) \ = \ \wpxioe 
\ee
solutions of which have {\it energy} which is independent of time $t$:
\be
\la{energy}
\sig\Bigl[\frac{\e\rho}{2} \dot{{\poseul}}_i^2 +\e \wxioe\Bigr] \ = \ E_0
\ee
where $E_0$ is determined by the initial data.  
The system \eqref{eom} has a {\it first order} in time formulation
obtained by setting:
\begin{align}
\label{fof}\begin{split}
  {\veleul}_i \ & = \ \dot {\poseul}_i\\
  \rho\frac{d{\veleul}_i}{dt} \ &= \ \frac{1}{\e} \wpxioe \,.
\end{split}
\end{align}

We are interested in studying the limit as $N\to \infty$, or
equivalently
$\e\to 0$, of this
system, and relating it to continuum elastodynamics. 
To this end we introduce by interpolation the following 
functions:
\begin{align}\la{int}
\begin{split}
y^\e (t,{\poslag}) & =  \sig \Big({\poseul}_i +\frac{1}{\e}({\poslag}- i\e))({\poseul}_{i+1} -{\poseul}_i)\Big)\id_{I^i_\e}({\poslag})\\
\tilde{y^\e} (t,{\poslag}) & = \sig {\poseul}_i\id_{I^i_\e}({\poslag})
\end{split}
\end{align}
for $I^{\e}_{i} = [i\e,(i+1)\e)$, as above.
We will prove that these two functions have the same limit as $\e\to
0$, and are thus lattice versions of the same macroscopic
object. In fact they are
lattice versions
of the Eulerian description of an elastic continuum, which
proceeds via a function ${\poslag}\mapsto
y(t,{\poslag})$ which gives the location in space of that infinitesimal 
part of the body whose
reference location is the point ${\poslag}$. It follows from the periodic
continuation  
${\poseul}_{N+i}(t) = {\poseul}_i(t)+2\pi \  \forall i$ that 
$y^\e(t,{\poslag}+2\pi)=y^\e(t,{\poslag})+2\pi$ and similarly for $\tilde{y^\e}$.
\begin{lemma}\la{uni}
Assume we have for each $N\in\{1,2,\dots\}$ a set of initial data
$\{({\poseul}_i(0),\dot {\poseul}_i(0))\}_{i=0}^{N-1}$ such that the energy is
uniformly bounded, so that \eqref{energy} with $\e=\frac{2\pi}{N}$
holds for some $E_0<\infty$ independent of $N$. Then for each such
$\epsilon$ the functions $y^\e$ and $\frac{\partial y^\e}{\partial t}$ are
bounded continuous functions of $t,{\poslag}$, and there exist a constant $C$
depending on the energy and on the coercivity constants, $C= C(E_0,\rho,
c_1, c_2)$ such that 
\begin{itemize}
\item[(i)] $\sup_t \Big( \|\frac{\partial y^\e}{\partial t} \|_{L^2
  (Q)}\ + \ 
\|\frac{\partial y^\e}{\partial {\poslag}} \|_{L^{p}(Q)}  \ +\ 
\|\frac{\partial \tilde{y}^\e}{\partial t} \|_{L^2 (Q)} \Big)
\ \leq \ C.$
\item[(ii)] $\sup_{t}\|\tilde y^\e - y^\e\|_{L^p (Q)} \leq C\e\,.$
\end{itemize}
\end{lemma}
\proof
Notice that  $|\frac{{\poslag} -i \e}{\e}\id_{I^i_\e}({\poslag})| \leq 1 $
everywhere. Therefore,
\bes
\| \frac{\partial y^\e}{\partial t} \|_{L^2(Q)}^{2}
  = \|\sig \Big(\dot{{\poseul}}_i + \frac{{\poslag}-i\e}{\e}(\dot{{\poseul}}_{i+1}
-\dot{{\poseul}}_i\Big)\id_{I^i_\e}\|_{L^2(Q)}^{2} \ \leq\  5\sig \e \|\dot{{\poseul}}_i\|^2_ {L^2(Q)}\\
  \ \leq \ C 
\ees
and similarly for $\tilde{y}^{\e}$. Next  
observe that $\frac{\partial y^\e}{\partial {\poslag}} 
=\sig\frac{{\poseul}_{i+1} -{\poseul}_i}{\e}\id_{I^i_\e}$ is bounded in $L^p$
by \eqref{energy} and our assumption on $W$, since 
$c_1\e|\frac{{\poseul}_{i+1}-{\poseul}_i}{\e}|^p\leq \e
W(\frac{{\poseul}_{i+1}-{\poseul}_i}{\e})+c_2$. This completes the proof of 
(i) using the energy bound \eqref{energy}.

The second assertion also follows from  \eqref{energy} and (i) since 
\bes
\tilde{y^\e} -y^\e \ =\  \sig \frac{({\poslag}_i -i\e)}{\e}
({\poseul}_{i+1}-{\poseul}_i)\id_{I^i_\e}\  = \ \sig \frac{({\poslag}_i -i \e)}{\e}\id_{I^i_\e} \frac{\partial y^\e}{\partial {\poslag}}\e
\ees
which implies (ii) as $|\frac{{\poslag}-i\e}{\e}\id_{I^i_\e}(\poslag)| \leq 1$.
\qed

For clarity, define the variables for the first order formulation,
\begin{align}\la{fov}
\begin{split}
u^\e(t,{\poslag}) &= \frac{\partial y^\e}{\partial {\poslag}}(t,{\poslag}) \\
v^{\e} (t,{\poslag}) &=\sig \dot{{\poseul}_i}\id_{I^i_\e} = \sig {\veleul}_i(t)\id_{I^i_\e} = \dtye\,.
\end{split}
\end{align}
Then the equations of motion \eqref{eom} in first order formulation 
\eqref{fof} become respectively,
\begin{align}
\begin{split}
   \e \rho \frac{\partial v^\e}{\partial t} & \ = \ 
 W'(u^\e(t,{\poslag})) - W'(u^\e(t, {\poslag}-\epsilon))\\
   \frac{\partial u^\e}{\partial t} & \  = \ \frac{\partial v^\e}{\partial
  {\poslag}}
-\frac{\partial}{\partial {\poslag}} (\dtye  - \dot {y}^\e)
\end{split}
\end{align}
which in weak form can be written as:
\begin{align}\la{wkf}\begin{split}
  &\int_0^\infty \!\int_{-\pi}^{\pi} \rho \frac{\partial \phi}{\partial t}v^\e - 
   \frac{1}{\e}\Bigl(\phi(t,{\poslag}+\e) -\phi(t,{\poslag})\Bigr) W'(u^\e ({\poslag})) \ d{\poslag}dt +
\int_{-\pi}^{\pi} \rho \phi({\poslag},0)v^\e({\poslag},0)\ d{\poslag}\ 
= 0\\
  & \int_0^\infty \!\int_{-\pi}^{\pi}  \Big(\frac{\partial
     \phi}{\partial t} u^\e - \frac{\partial\phi}{\partial {\poslag}} v^\e
-\frac{\partial^2
     \phi}{\partial {\poslag}\partial t} ({\tilde{y}}^\e  - {y}^\e)
\Big) \ d{\poslag}dt \ + \
\int_{-\pi}^{\pi}\phi({\poslag},0)u^\e({\poslag},0)\ d{\poslag} \
= 0 
\end{split}
\end{align}
for all $\phi\in C_c^2(Q_\infty)$.

As in lemma \ref{uni} bounds which are uniform in $\e$ come 
from energy conservation, which in first order variables takes the form
\be
\frac{\rho}{2} \|v^\e(t,\,\cdot\,)\|^{2}_{L^2} + 
\int W(u^\e(t,\,\cdot\,)) d{\poslag} \ = \
\frac{\rho}{2} \|v^\e(0,\,\cdot\,)\|^{2}_{L^2} + 
\int W(u^\e(0,\,\cdot\,)) d{\poslag} \ 
\leq E_0 <\infty.
\ee
Thus
$
\sup_t\bigl(\int |u^\e|^p d{\poslag} +
\int |v^\e|^2 d{\poslag}\bigr) \leq  C\,.
$
To take the limit of \eqref{wkf} we use the facts that
$
\frac{\phi(t,{\poslag}+\e) -\phi(t,{\poslag})}{\e}  \longrightarrow \frac{\partial
  \phi}{\partial {\poslag}}
$ uniformly (since $\phi$ is a test function) and 
$
\tilde{y^\e} -y^\e \longrightarrow 0 
$
in $L^p$ by lemma \ref{uni}.

In the limit $\e\to 0$ there is a Young measure $\n$ which 
represents
weak limits of the sequence $(u^\e,v^\e)\wkarr (u,v)$:
$$
v=\int\,\lambda\,d\n(M,\lambda)\quad\hbox{and}\quad
u=\int\,M\,d\n(M,\lambda)\,,
$$
and 
of functions $g(u^\e,v^\e)$ which are $L^1$ precompact, so that
in particular
\begin{align}
\lim_{\e\to 0}\int_0^\infty\!\int_{-\pi}^{\pi} \phi g(u^\e,v^\e)\,d{\poslag}dt
&=
\int_0^\infty\!\int_{-\pi}^{\pi} \phi\,\langle\n,g\rangle\,d{\poslag}dt\\
&=
\int_0^\infty\!\int_{-\pi}^{\pi} \int\, \phi\, 
g(M,\lambda)\,d\n(M,\lambda)\,d{\poslag}dt
\end{align}
for all bounded $\phi$ which are 
$2\pi$-periodic in ${\poslag}$ and vanish for large $t$. On the other hand
for the energy density $\eta(u^\e,v^\e)=\frac{\rho}{2}(v^\e)^2+W(u^\e)$
we only have $L^1$ boundedness, and the weak limit includes a concentration
measure $\g$:
$$\lim_{\e\to 0}
\int_0^\infty\!\int_{-\pi}^{\pi} \phi \eta(u^\e,v^\e)\,d{\poslag}dt
=\int_0^\infty\!\int_{-\pi}^{\pi} \phi\, \langle\n,\eta\rangle\,d{\poslag}dt
+\int_0^\infty\!\int_{-\pi}^{\pi} \phi\,\g(d{\poslag}dt)
$$
for $\phi\in C_c(Q_\infty)$.
Consider initial data with the property that 
$(u^\e({\poslag},0),v^\e({\poslag},0))\,\to\,(u({\poslag},0),v({\poslag},0))$
in $L^p\times L^2$, and $\int\,\eta(u^\e({\poslag},0),v^\e({\poslag},0))\,dX\,
\to\,\int\,\eta(u({\poslag},0),v({\poslag},0))\,dX$.
On account of the
assumptions on $W$ the limit $(u,v,\n)$ is a dissipative
measure-valued solution in the sense that:
\begin{align}
\int_0^\infty\!\int_{-\pi}^{\pi} 
              \Big(\rho v\,{\partial_t \phi}
                   +\langle\n, W'\rangle \,{\partial_{\poslag} \phi}\Big) 
                     \ d{\poslag}dt \ & 
 + \
\int_{-\pi}^{\pi} \rho \phi({\poslag},0)v({\poslag},0)\ d{\poslag}\ 
= 0 \\
\int_0^\infty\!\int_{-\pi}^{\pi} 
              \Big( u\,{\partial_t \phi} 
                   -v\,{\partial_{\poslag} \phi}\Big) \ d{\poslag}dt \ &
+ \
\int_{-\pi}^{\pi}\phi({\poslag},0)u({\poslag},0)\ d{\poslag} \ = 0 \,,
\end{align}
for all $\phi\in C_c^1(Q_\infty)$, and 
\beq
\la{decont}
\int_0^\infty\!\int_{-\pi}^{\pi} \dot \theta(t)
\bigl(\langle\n,\eta\rangle\,d{\poslag}dt
+\g(dXdt)\bigr)
+\theta(0)\!\int_{-\pi}^{\pi} \eta(\ou(\poslag,0)\,,\ov(\poslag,0)) \,d{\poslag}
\geq 0\,,
\eeq
for non-negative $\theta\in C^1([0,\infty))$. 
(In fact the dissipative condition 
\eqref{decont} holds as an equality.)

Now using the relative entropy method and the convexity assumption 
on $W$ we can prove that in fact the convergence is strong and
concentration free when
a classical solution $(\ou,\ov)$ exists on $\barQT$:
\begin{theorem}
Assume that there is a pair of Lipschitz functions 
$(\ou,\ov)\in W^{1,\infty}(\barQT)$ which satisfy the continuum limit equations:
\begin{align}\la{wkfe}\begin{split}
  &\int_0^\infty \!\int_{-\pi}^{\pi} \rho \frac{\partial \phi}{\partial t}\ov + 
   \frac{\partial \phi}{\partial {\poslag}} W'(\ou ({\poslag})) \ d{\poslag}dt 
 +
\int_{-\pi}^{\pi} \rho \phi({\poslag},0)\ov({\poslag},0)\ d{\poslag}\ 
= 0\\
  & \int_0^\infty \!\int_{-\pi}^{\pi}  \Big(\frac{\partial
     \phi}{\partial t} \ou - \frac{\partial\phi}{\partial {\poslag}} \ov
\Big) \ d{\poslag}dt 
+ \
\int_{-\pi}^{\pi}\phi({\poslag},0)\ou({\poslag},0)\ d{\poslag} \
= 0 \,,\end{split}
\end{align}
for all $\phi\in C_c^1(Q_T)$.
Assume that there is a sequence of initial configurations of the
lattice $\{({\poseul}_i(0),\dot {\poseul}_i(0))\}_{i=0}^{N-1}$ with uniformly
bounded energy, and such that the
corresponding interpolated functions 
$(u^\e({\poslag},0),v^\e({\poslag},0))$, $\e=\frac{2\pi}{N}$, 
converge strongly to
$(\ou({\poslag},0),\ov({\poslag},0))$
in $L^p\times L^2$ and 
$\int\,\eta(u^\e({\poslag},0),v^\e({\poslag},0))\,dX\,
\to\,\int\,\eta(\ou({\poslag},0),\ov({\poslag},0))\,dX$\,.
Then $(u^\e,v^\e)$, as defined in \eqref{int} and 
\eqref{fov} from the
solutions $\{({\poseul}_i(t),\dot {\poseul}_i(t))\}_{i=0}^{N-1}$ of the 
microscopic model, converge strongly in $L^p\times L^2(Q_T)$ 
to the continuum limit $(\ou,\ov)$. 
Alternatively said, the Young measure $\n$
is a Dirac measure supported on
$(\ou,\ov)$ and there is no
concentration, i.e. the concentration measure $\g$ is null.
\end{theorem}
\proof
We define the relative entropy as 
$h(\n,u,v,\ou,\ov)=\langle\n,\eta_{rel}\rangle
=\int\eta_{rel}(M,\lambda;\ou,\ov)\,d\n(M,\lambda)$
with 
\begin{align}
\eta_{rel}(M,\lambda;\ou,\ov)&=\eta(M,\lambda)-\eta(\ou,\ov)-\ov(\lambda-\ov)
-W'(\ou)(M-\ou)\notag\\
&=\frac{\rho}{2}(\lambda-\ov)^2+W(M)-W(\ou)-W'(\ou)(M-\ou)
\,.\notag
\end{align}
Under the assumptions on $W$ above there exists $C>0$ such that
$$
\frac{W'(M)-W'(\ou)-W''(\ou)(M-\ou)}
{W(M)-W(\ou)-W'(\ou)(M-\ou)}\leq C
$$
everywhere. (The number $C$ depends upon the bounded region $D$ in which
$\ou$ takes its values).
Given this inequality and the assumption that the intial data converge to the 
initial data $(\ou_0,\ov_0)$ of the bounded 
Lipschitz  solution $(\ou,\ov)$ we then deduce, via a
calculation  analogous to that in \eqref{dm}-\eqref{weceldr1}, that
$$
\int\,h(\n,u,v,\ou,\ov)\,d{\poslag}\bigr|_{t}\leq C'\int_0^t\int\,h(\n,u,v,\ou,\ov)\,d{\poslag}
\bigr|_{\tau}d\tau
$$
for $0\leq t<T$, and hence that $h$ and $\g$ are zero almost 
everywhere for positive times for which the classical solution exists. This
implies that $\n_{({\poslag},t)}  = \delta_{(u({\poslag},t),v({\poslag},t))}$ as 
previously, and hence that the convergence of $(u^\e,v^\e)$ to $(u,v)$ 
is strong and concentration free as claimed.
\qed

\setcounter{section}{0}
\setcounter{equation}{0}
\protect\renewcommand{\thesection}{\Alph{section}}

\section{Appendix: An energy concentration measure for 
measure-valued solutions}\label{concen}

In this appendix we summarize what we need about 
the Young measure description of
oscillations and concentrations in weakly convergent sequences
of functions $f^\epsilon(y)\in\R^m$ defined on the set $\barQT$, writing 
$y$ as the independent variable ($y=(t,x)$). 

We consider two settings in which the Fundamental Therem of Young Measures,
as found in  Ball \cite{ball88}, applies: the $L^\infty$ setting of section
\ref{linf} and the  $L^p$ setting of sections \ref{qel} and \ref{pel}.
In the $L^\infty$ setting the theorem attaches to a uniformly bounded
sequence of functions on $\barQT$ a subsequence, still written $f^\e$, and a
parametrized Young measure (meaning a weak* measurable
$\barQT$-parametrized family of probability measures $\n=(\n_y)_{y\in\barQT} $)  
such that
for any continuous function $F: \R^m \to \R$  
\beq\la{ymrep}
F(f^\e) \wkarr \langle \n\,, F\rangle \qquad\mbox{weak*\ in \ } L^\infty (\barQT). 
\eeq
In the $L^p$ setting, $1<p<\infty$, a similar conclusion holds for any
sequence of functions $f^\epsilon$ which are bounded in $L^p$: 
for  
continuous $F$
such that $F(f^\e)$ is $L^1$ weakly precompact there holds
\beq\la{ymprep}
F(f^\e) \wkarr \langle \n\,, F\rangle \qquad\mbox{weakly\ in \ } L^1 (\barQT). 
\eeq
This representation will generally not hold if {\em $L^1$ weakly
precompact}
is replaced by {\em $L^1$ bounded} because concentrations can develop.
Various tools have been introduced to
describe this such as biting convergence,
the generalized concentration Young measure, microlocal defect
measure, $H$-measure, varifold measure included, see references
\cite{dipm, ballmurat, ger89,tartar90,ab97,fmp98} and \cite[Section 1.D]{evans}.
Here we introduce by hand a simple  measure 
$\g$ of concentration effects {\it in the
energy} or other non-negative functions $F$ of {\em critical growth}, that
is functions such that $F(f^\e)$ is bounded, but not necessarily 
weakly precompact, in $L^1$
(for example,  $|f^\e|^p$  of an $L^p$-bounded sequence).   This measure
$\g$ is a sharpening of the weak* defect measure $\s$ of Lions 
(described in  (\cite[Chapter 9]{mb}). In fact its existence follows 
as a particular case of a quite general result \cite[Theorem 2.5]{ab97}. 
However since
we only need a rather special case - to describe the weak limit of
a single non-negative function $\eta$ - we give a simple direct proof from 
first principles.

We introduce this measure in two separate cases, first 
for illustrative purposes in the $L^2$ setting 
which applies in section \ref{qel}, and then in the more
general setting which is useful in the case of a polyconvex  energy of
section \ref{pel}.


\subsection{The $L^2$ case}\la{el2}
We now consider the case $p=2$ in more detail:
let $f^\epsilon(y)$ converge weakly in $L^2$ to $f(y)$, and assume that
$\int|f^\epsilon(y)|^2dy\le K<\infty$.
Then by the previous discussion
$$
\int_{Q_T}  F(f^\e)(y) w(y) dy \lra \int_{Q_T} \langle {\n}_y , F\rangle w(y)dy 
$$
for any $w\in L^\infty (\barQT)$ and for any $F$ satisfying 
$\lim_{|z|\to +\infty}\frac{|F(z)|}{1+|z|^2}=0$, (since this implies that
$F(f^\e)$ is
sequentially weakly precompact in $L^1$ by the criterion of   de la
Vallee Poussin.)
For the function $F(z)=|z|^2$ itself, however,  
$y\mapsto |f^\epsilon(y)|^2dy$ are
weak* precompact in the space of non-negative Radon measures $\cM^+(\barQT)$,
and the functions $|f^\epsilon|^2$ need not be weakly precompact in $L^1$ and 
as a result the Young measure representation in general fails. 

  In this context we define a defect measure   
by applying the Banach-Alaoglu theorem to the sequence
$|f^\epsilon-f|^2$ to obtain a subsequential weak* limit $\s$,
which is a non-negative Radon measure,
\beq
\s(\psi)=\iint\psi d\s=\lim_{\epsilon\to 0}
\iint\,\psi\,|f^\epsilon-f|^2
\,dxdt\,,
\eeq
for all $\psi\in C(\barQT)$. Alternatively, noting the identity
$|f^\epsilon|^2=|f^\epsilon-f|^2+|f|^2+2\langle f, f^\epsilon-f\rangle$,
it follows from the definition of weak $L^2$ convergence that an equivalent
definition is $$\s= \hbox{wk*-}\lim_{\epsilon\to 0}(|f^\epsilon|^2-|f|^2)\in
\cM^+(\barQT).$$ Simple examples indicate that $\s$ can be non-zero due
to purely oscillatory effects, and it is ``too large'' to describe
concentration effects in a useful way. Therefore
we will use a modification of $\s$,
called $\g$, which is smaller (i.e. $\g(E)\leq\s(E)$) and is designed
to be useful to describe weak limits of non-negative functions of
critical growth. 
To introduce the measure $\g$ we
first observe that if we apply the Young measure theorem to $f^\epsilon$
we obtain for almost every $y\in \barQT$ a Radon probability measure
$\n_{y}$, and the function $\int|\lambda|^2\n_{y}(d\lambda)$ is well defined
in the extended non-negatives $[0,\infty]$ by the
monotone convergence theorem. Indeed let 
$q_R(\lambda)=|\lambda|^2\id_{|\lambda\leq R}
+R^2\id_{|\lambda|\geq R}$ then $q_R(\lambda)\nearrow q(\lambda)=|\lambda|^2$
and so $\int|\lambda|^2\n_{y}(d\lambda)=\lim_{R\to\infty}
\int q_R(\lambda)\n_{y}(d\lambda)
$
is well defined for a.e. $y$ and is in $L^1(\barQT)$ 
since by the Young measure representation
theorem
$$
\iint \psi(y)q_R(\lambda)\n_{y}(d\lambda)dy
=\lim\int\psi(y)q_R(f^\epsilon(y))dy
\leq K\max_{y\in Q}|\psi(y)|
$$
for all $\psi\in C(\barQT)$; choosing $\psi(y)\equiv 1$ allows us to apply the
monotone convergence theorem again to deduce that 
$\langle\n_y(\lambda), |\lambda|^2\rangle=
\int |\lambda|^2\n_{y}(d\lambda)\in L^1(\barQT)$ since it is a monotone
non-decreasing
limit of non-negative functions of uniformly bounded integral.

Now to define the concentration measure $\g$, we just mimic the definition
of $\sigma$, replacing $|f(y)|^2$ by 
$\langle\n_y(\lambda), |\lambda|^2\rangle$, i.e. we consider 
$\hbox{wk*-}\lim_{\epsilon\to 0}(|f^\epsilon(y)|^2-\langle\n_y(\lambda), 
|\lambda|^2\rangle)$. To show that this limit exists in $\cM^+(\barQT)$ we
use again the Young measure representation: for any $R>0$ and any 
{\em non-negative} function $\psi\in C(\barQT)$,
$$
\iint \psi(y)q_R(\lambda)\n_{y}(d\lambda)\,dy=
\lim_{\epsilon\to 0}\int\psi(y)q_R(f^\epsilon(y))\,dy
\leq \lim_{\epsilon\to 0}\int\psi(y)|f^\epsilon(y)|^2\,dy
$$
and therefore
$$
\iint \psi(y)|\lambda|^2\n_{y}(d\lambda)\,dy=
\sup_{R>0}\iint \psi(y)q_R(\lambda)\n_{y}(d\lambda)\,dy
\leq \lim_{\epsilon\to 0}\int\psi(y)|f^\epsilon(y)|^2\,dy
$$ and hence \beq\la{dg2} \g= \hbox{wk*-}\lim_{\epsilon\to
  0}\bigl(|f^\epsilon|^2-\langle\n_y(\lambda),
|\lambda|^2\rangle\bigr) \in\cM^+(Q) \eeq is a well defined {\em
  non-negative} Radon measure. Since H\"older's inequality implies
that $|f(y)|^2=|\langle\n_y,\lambda\rangle|^2\leq
\langle\n_y,|\lambda|^2\rangle$, this definition implies that
$\g\leq\s$ as claimed earlier. The reason that the concentration Young
measure $\g$ is useful is that it allows a description of the weak
limit of the energy, in terms of the Young measure $\n$ - the measure
defined in \eqref{dg2} is used in section \ref{qel}.

\subsection{The general case}\la{elp}
To describe concentration
effects arising from more general energy functionals, such as the
polyconvex ones in section \ref{pel}, it is necessary to generalize
the preceding definition.  We now show that the same argument can be
applied to any {\em non-negative} continuous function $\eta$ which satisfies
$\int \eta(f^\e)\le K<\infty$, but for which the de la Vallee Poussin criterion
does not apply and weak $L^1$ precompactness of $\eta(f^\e)$ cannot be
assumed.  Instead we assume that $\eta\geq 0$ is a superlinear
function and $\sup_{\e >0}\int
\eta(f^\e) dx < K$ where $f^\e $ is assumed to be a sequence of
Lebesgue measurable functions which according to the theorem of Ball
(\cite{ball88}) has a subsequence, also called $f^\e$, with
associated Young measure $\n_y$, which is a weak* measurable family of
Radon {\em probability} measures on account of the superlinearity assumption
on $\eta$.  By the same theorem the Young measure represents $L^1$ weak
limits of compositions of the $f^\e$ as in \eqref{ymprep} when these are $L^1$ precompact.
Observe that 
$$
y \mapsto \int \eta(\lambda)\n_{y}(d\lambda)
$$ 
is well defined a.e. in $y$ with values  in the
extended non-negatives $[0,\infty]$ and is in $L^1$ by the monotone convergence
theorem:
$\eta_R(\lambda)=\eta(\lambda)\id_{\eta(\lambda)\leq R}
+R\id_{\eta(\lambda)\geq R}$ then $\eta_R(\lambda)\nearrow
\eta(\lambda)$ and so $\int
\eta(\lambda)\n_{y}(d\lambda)=\lim_{R\to\infty} \int
\eta_R(\lambda)\n_{y}(d\lambda) $ is well defined for a.e. $y$ and is
in $L^1(\barQT)$ since by the Young measure representation theorem
$$
\iint \psi(y)\eta_R(\lambda)\n_{y}(d\lambda)dy
=\lim_{\e \to 0}\int\psi(y)\eta_R(f^\epsilon(y))dy
\leq K\max_{y\in \barQT }|\psi(y)|
$$
for all $\psi\in C(\barQT)$. Choosing $\psi(y)\equiv 1$ allows us to
apply the monotone convergence theorem again to deduce that   
\be\la{etalone}
\langle\n_y(\lambda), \eta(\lambda)\rangle=
\int \eta(\lambda))\n_{y}(d\lambda)\in L^1(\barQT)
\ee
since it is a non-decreasing limit of non-negative functions of uniformly bounded integral: explicitly, by the Young measure representation for 
$\eta_R (f^\e)$,
$$
\int \langle\n_y (\lambda) , \eta_R (\lambda) \rangle = \lim_{\e\to 0}\int\eta_R(f^\e) \ \leq \ \sup_{\e} \int\eta_R (f^\e) \ \leq \sup_{\e} \int \eta(f^\e) \leq K 
$$
by assumption on $\eta$ and $(f^\e)$ where the integrals are over
$Q_T$ and using that  $0<\eta_R \nearrow \eta$
we deduce \eqref{etalone} by monotone convergence taking the limit in $R$ . 

It is not, however, the case that $\eta(f^\e) $ are $L^1$ precompact and
so $\langle \n_y, \eta\rangle$ does not give its weak limit in general 
due to concentration.
The concentration effect can be measured by considering 
the {\em concentration measure}
 \be \la{conc}
\g = \hbox{wk*-} \lim_{\epsilon \to 0} (\eta(f^\e)-\langle \n_y
,\eta\rangle) \ee which is a well defined non-negative Radon measure for a subsequence of the $\eta(f^\e)$ (since they have bounded integral):
to see that $\g$ is indeed non-negative we
use again the Young measure representation to deduce that for any
$R>0$, and any {\em non-negative} function $\psi\in C(\barQT)$,
$$
\iint \psi(y)\eta_R(\lambda)\n_{y}(d\lambda)\,dy=
\lim_{\epsilon\to 0}\int\psi(y)\eta_R(f^\epsilon(y))\,dy
\leq \lim_{\epsilon\to 0}\int\psi(y)\eta(f^\epsilon(y))\,dy
$$
and therefore
$$
\iint \psi(y)\eta(\lambda)\n_{y}(d\lambda)\,dy=
\sup_{R>0}\iint \psi(y)\eta_R(\lambda)\n_{y}(d\lambda)\,dy
\leq \lim_{\epsilon\to 0}\int\psi(y)\eta(f^\epsilon(y))\,dy
$$
and hence
$$
\g=
\hbox{wk*-}\lim_{\epsilon\to 0}\bigl(\eta(f^\epsilon )-\langle\n_y, 
\eta\rangle\bigr)  \in\cM^+(\barQT)
$$
is a well defined {\em non-negative} Radon measure. 

If in addition $\eta$ is convex,
then $\g\leq\s$ where $\s$ is the natural generalization of the weak*
defect measure, namely $\s=\hbox{wk*-}\lim_{\epsilon\to
  0}(\eta(f^\epsilon)-\eta(f))\in \cM^+(\barQT).$  This is an immediate
consequence of  Jensen's
inequality which  implies that 
$$
\int \eta(\lambda)\ d\n \geq \eta(\int
\lambda\ d\n)=\eta( \lim f^\e)= \eta(f).
$$

 The reason that
$\g$ is useful is that it allows a description of the weak limit of
the energy, in terms of the Young measure $\n$.  In section \ref{pel}
 this applies  to a sequence $f^\e = (v^\e, \Xi^\e) $ which is 
 bounded in a direct sum of
 different Lebesgue spaces, and which therefore has a weak limit point
 in the same space.
\begin{remark}
Although we refer to $\g$ as concentration measure, it is not always
supported on a small set: there exist sequences of functions in which
the concentration smears out to fill the whole domain, see
\cite[Example 2]{ballmurat}.
\end{remark}

\noindent
{\bf Acknowledgements} 
This research was started while AET was visiting the Newton Institute, and
completed with support by the EU FP7-REGPOT project "Archimedes Center for
Modeling, Analysis and Computation", during a visit
by DS and SD. Also supported by the National Science Foundation
and EPSRC.

\end{document}